\documentclass[12pt]{article}
\usepackage{color}   
\usepackage{amsmath, amssymb}
\usepackage[latin2]{inputenc}
\usepackage{t1enc}
\usepackage{bbm}
\usepackage{auto-pst-pdf}
\usepackage{pstricks}
\usepackage{pstricks-add}
\newtheorem{lem}{Lemma}

\newtheorem{prop}{Proposition}
\newtheorem{thm}{Theorem}
\newtheorem{cor}{Corollary}

\newcommand{\bee}{\begin{equation}}
\newcommand{\eee}{\end{equation}}
\newcommand{\bees}{\begin{equation*}}
\newcommand{\eees}{\end{equation*}}
\newcommand{\bali}{\begin{aligned}}
\newcommand{\eali}{\end{aligned}}

\newcommand{\dss}{\:\mathrm{d}\mathbf{s}}
\newcommand{\dxx}{\:\mathrm{d}\mathbf{x}}
\newcommand{\dyy}{\:\mathrm{d}\mathbf{y}}
\newcommand{\dzz}{\:\mathrm{d}\mathbf{z}}
\newcommand{\dx}{\:\mathrm{d}x}

\newcommand{\dr}{\:\mathrm{d}r}
\newcommand{\pa}{\partial}
\newcommand{\er}{\mathbb{R}}

\newcommand{\nnu}{\textrm{\boldmath$\nu$}}
\newcommand{\pphi}{\textrm{\boldmath$\phi$}}
\newcommand{\xxi}{\textrm{\boldmath$\xi$}}

\newcommand{\IP}{\textrm{IP}}
\newcommand{\BV}{\textrm{BV}}

\newcommand{\lm}{\lesssim}
\newcommand{\gm}{\gtrsim}

\newcommand{\nnull}{\mathbf{0}}
\newcommand{\mcd}{\mathcal{D}}
\newcommand{\mbp}{\mathbb{P}}
\newcommand{\mcf}{\mathcal{F}}
\newcommand{\mct}{\mathcal{T}}
\newcommand{\mck}{\mathcal{K}}

\newcommand{\kk}{\textbf{k}}
\newcommand{\ww}{\mathbf{w}}
\newcommand{\xx}{\mathbf{x}}
\newcommand{\yy}{\mathbf{y}}
\newcommand{\zz}{\mathbf{z}}
\newcommand{\hh}{\mathbf{h}}
\newcommand{\sss}{\mathbf{s}}
\newcommand{\ordo}{\mathcal{O}}

\newcommand{\jump}[1]{\ensuremath{\left[ \! \left[ #1 \right] \! \right]} }
\newcommand{\aver}[1]{\ensuremath{\left\{ \!\! \left\{ #1 \right\} \!\! \right\}}}

\newcommand{\proof}{\emph{Proof:}\:}
\newcommand{\eproof}{\quad$\square$}

\title{Energy norm error estimates for averaged discontinuous Galerkin 
methods: multidimensional case}

\author{Ferenc Izs\'ak
\thanks{Department of Applied Analysis and Computational 
 Mathematics, Faculty of Science, E\"otv\"os Lor\'and University,
P\'azm\'any s\'et\'any 1/C, H-1117 Budapest, Hungary. 
Email: \texttt{izsakf@cs.elte.hu} Tel.:+3613722500/8428}} 

\begin{document}
\maketitle
\begin{abstract}
A mathematical analysis is presented for a class of interior penalty (IP)
discontinuous Galerkin approximations of elliptic boundary value problems.
In the framework of the present theory one can derive some overpenalized 
IP bilinear forms in a natural way avoiding any heuristic choice of fluxes 
and penalty terms.
The main idea is to start from bilinear forms for the local average of 
discontinuous approximations which are rewritten using the theory
of distributions. It is pointed out that a class of overpenalized 
IP bilinear forms can be obtained using a lower order perturbation of these. 
Also, error estimations can be derived between the local averages of
the discontinuous approximations and the analytic solution in the 
$H^1$-seminorm. Using the local averages, the analysis is performed in
a conforming framework without any assumption on extra smoothness 
for the solution of the original boundary value problem.  
\end{abstract}

\section{Introduction}
Discontinuous Galerkin (dG) methods have been introduced and used from
the end of the seventies, first for linear transport problems. Later 
this was generalized to elliptic boundary value problems and nowadays 
it is available for the numerical solution of almost all kind 
of problems based on PDE's.  

These methods have proved their usefulness in several simulations of 
real-life phenomena \cite{dawson11}, \cite{hartmann10}, \cite{tago12}. 
The most favorable property of the corresponding 
numerical methods is that the local mesh refinement can easily be 
performed giving rise to efficient adaptive strategies.

An important milestone in the systematic analysis of dG methods for the
elliptic boundary value problems was the paper \cite{arnold01}. This 
pioneering work served as a basis of the consecutive works concerning a 
priori and a posteriori error estimates \cite{karakashian03}, 
$hp$-adaptive methods \cite{houston07b}, time dependent problems 
\cite{cockburn98}. For an up-to-date summary of the theoretical achievements 
for dG methods we refer the recent monograph \cite{dipietro12} and for
 implementation issues the monographs \cite{hesthaven07} and \cite{riviere08}.

At the same time, the above analysis should be improved in some aspects. 
First, which can be considered as a didactic issue, the choice of the 
corresponding bilinear forms would deserve more motivation. After recasting 
the elliptic problem in a mixed form, numerical fluxes and penalty terms
are defined which lead to different bilinear forms. No a priori
suggestion or motivation (on a physical basis) is mentioned to 
propose an appropriate choice of the fluxes. 
A similar situation arises when penalty terms are defined. 

The second issue is the assumption on extra-regularity of the analytic
solution. This problem was solved in the meantime: in \cite{gudi10}
the author developed an analysis based on a Strang type lemma 
\cite{ern04}, which could successfully deal with the 
non-conformity of the dG type approximation.
   
The most important issue is the norm for the convergence. The
choice of the bilinear form implies a mesh-dependent norm, which is a
real mathematical artifact. The convergence is proved with respect to 
this norm or in a weaker, \emph{e.g.}, in the $L_2$-norm. At the same 
time, in the corresponding real-life problems the natural norm is usually
the $H^1$-norm (or seminorm). Note that there are some achievements 
which point out the usefulness of the interior penalty (IP) methods. 
For these methods, one can obtain convergence in the so-called BV
 norm which does not depend on the actual mesh \cite{buffa09},
\cite{dipietro10} and can be related to broken Sobolev norms.

The aim of the present work is to contribute to the mathematics of the
dG methods for elliptic boundary value problems by proposing an alternative 
of the commonly used theoretical basis in \cite{arnold01}. 
In particular, we derive overpenalized interior penalty bilinear forms
in a natural way avoiding the notion of numerical fluxes or recasting 
them into a mixed form.
The new idea is to use the local average of the discontinuous 
approximation from the beginning. 
The main benefit of the analysis is that it can be done in an $H^1$-conforming 
framework such that one can prove the quasi optimal convergence of the 
local average with respect to 
the natural $H^1$-seminorm for Dirichlet problems. 
This work is a generalization of the paper \cite{csorgo14} concerning 
the one-dimensional case. 

The idea to use \emph{post}processing (or smoothing or filtering) for dG 
approximations has already appeared in the literature \cite{cockburn03b}. 
In the last years, many related results have been achieved: 
involved algorithms were developed for linear hyperbolic problems 
in \cite{king12} and their 
accuracy-increasing property was verified also for advection-diffusion 
problems with respect to negative Sobolev norms \cite{ji12}. 
The accurate computation of the corresponding convolutions is challenging,
see the recent developments in \cite{mirzaee13a} and \cite{mirzaee14}.  


The setup of the article is as follows.
After some preliminaries we give the bilinear form for the averaged 
approximation, which still contains convolution terms. 
We then expand the terms and point out that with a lower-order 
perturbation an overpenalized IP bilinear form can be obtained. 
This result is given in Theorem 1. 
Based on this, we can state the closedness of the approximation from 
the new bilinear form and the one arising from the overpenalized IP
bilinear form, see Theorem 2. Finally, in Theorem 3, an optimal convergence 
rate for the averaged overpenalized IP approximation is proved in the $H^1$ 
(semi)norm. 
The only tool we use beyond the standard armory of the finite element 
analysis is a bit of distribution theory.

\section{Mathematical preliminaries}
We investigate the finite element solution of the elliptic boundary
value problem
\bee\label{basic_eq}
\begin{cases}
-\Delta u (\xx) = g(\xx)\quad \xx\in\Omega\subset\er^d\\
u(\xx) = 0\quad \xx\in\pa\Omega,
\end{cases}
\eee
where $\Omega$ is a polyhedral Lipschitz-domain and $g\in L_2(\Omega)$
is given.

The finite element approximation is computed on a non-degenerated 
simplicial mesh $\mct_h $ with the mesh parameter $h$. The symbol
$\mcf$ denotes the set of interelement faces.
For the numerical solution we use the finite element space 
\[
\mbp_{h,\kk} = \{u\in L_2(\Omega): u|_{K}\in P_{k_j}(\Omega_j)\;
\textrm{for all} \;\Omega_j\in\mct_h\},
\] 
where $\kk = (k_1,k_2,\dots)$ and $\mathbb{P}_{k_j}(\Omega_j)$ denotes the 
linear space of polynomials of total degree $k_j$ on the subdomain 
$\Omega_j$. This notation will also be used for interelement faces 
and for balls instead of $\Omega_j$. 
We also make use of the conventional notation 
$\aver{\cdot}:\mbp_{h,\kk}\to L_2(\mcf)$ and 
$\jump{\cdot}:\mbp_{h,\kk}\to L_2(\mcf)$ for the average and jump 
operators which are given on each interelement face 
$f_\Omega = \bar\Omega_+\cap\bar\Omega_-$ with 
$$
\aver{v}_{f_\Omega} (\xx) = \frac{1}{2} (v(\xx_+) + v(\xx_-))
\quad\textrm{and}\quad
\jump{v}_{f_\Omega} (\xx) = \nnu_+ v(\xx_+) + \nnu_- v(\xx_-).
$$ 
Here $\nnu_\pm$ denotes the outward normal of $\Omega_\pm$ and 
$v(\xx_\pm) = \lim_{\Omega_\pm\supset\xx _n\to\xx} v(\xx _n)$. 
On each boundary face $f\subset \pa\Omega$ we simply define 
$$
\aver{v}_{f} (\xx) = v(\xx)
\quad\textrm{and}\quad
\jump{v}_{f} (\xx) = \nnu(\xx) v(\xx).
$$ 
The $L_2(\Omega^*)$ norm on a generic domain $\Omega^*$ will be denoted 
with $\|\cdot\|_{\Omega^*}$ and the corresponding scalar product with 
$(\cdot, \cdot)_{\Omega^*}$. In case of $\Omega^* = \Omega$ or if the 
support of the terms is given, we omit the subscript. Similar notation 
is applied for the scalar product and the corresponding  $L_2$ norm on 
$\mcf$ and on a single interelement face $f$. 

With these, the most popular dG approximation of $u$ in \eqref{basic_eq}
is the so-called symmetric interior penalty dG method which is given with 
the bilinear form 
$a_{\textrm{IP}}: \mathbb{P}_{\hh,k}\times \mathbb{P}_{\hh,k}\to \er$ as 
follows:
\bee\label{IP_form}
a_{\textrm{IP}}(u,v) = 
(\nabla_h u, \nabla_h v) 
-\sum_{f\in\mathcal{F}}
(\aver{\nabla_h u}, \jump {v})_f + (\aver{\nabla_h v}, \jump {u})_f
+\sum_{f\in\mathcal{F}} \sigma_h (\jump{u}, \jump {v})_f,
\eee
where $\nabla_h$ denotes the piecewise gradient on the subdomains in 
$\mathcal{T}_h$ and $\sigma_h\in\er$ denotes a penalty parameter,
which is proportional with $(\textrm{diam}\:f)^{-1}$ in the conventional 
setting.  We will also use the notation $\nabla_f \jump{u}$ for 
the gradient of the jump functions defined on the interelement face $f$. 

The notation $\lambda_d(\cdot)$ will be used to the $d$-dimensional
Lebesque measure. For the local average we use the piecewise constant 
function $\eta_h: \er^n\to\er$ depending also on the parameter $s>1$ with
$$
\eta_h (\xx) =
\begin{cases}
\frac{1}{B_{h^s,d}}\quad |\xx|\le h^s\\
  0\quad |\xx|> h^s,
\end{cases}
$$
where $B(\xx, r)$ denotes the closed ball with radius $r$ centered at 
$\xx$ and $B_{h^s,d} = \lambda_d(B(\nnull, h^s))$. The analysis makes use 
only two properties of $\eta_h$: this is symmetric with respect to the 
origin and $\int_{\er^d}\eta_h = 1$ such that $\eta_h * u$ is the local 
average of the function $u:\er^d\to\er$.  
Also, a straightforward computation gives that
$\textrm{supp}\:\eta_h * \eta_h = B(\nnull, 2h^s)$ and
$\int_{B(\nnull, 2h^s)} \eta_h * \eta_h = 1$. These facts will be 
used without further reference.

The analysis of the conforming approach will be carried out in the space 
$$
\mbp_{h,\kk,s} = 
\{\eta_h * u_0|_{\Omega_h}: u_0\; \textrm{is the zero extension of}\; u \in \mbp_{h,\kk}\},
$$
where $\Omega_h = \{\xx\in\er^d: d(\xx,\Omega)< h^s\}$.
Obviously, $\mbp_{h,\kk,s}\subset H_0^1(\Omega_h)$. We use the notation 
$\Omega_{j,h}$ in a similar sense and $\tilde \Omega_j = 
\textrm{int}\:\{ \bar\Omega_k\in \mathcal{T}_h: 
\bar\Omega_j\cap\bar\Omega_k \not= \emptyset\}$ 
for the patch of $\Omega_j$.

\begin{figure}
\begin{pspicture}(0,-1)(10, 6)

\psarc[linewidth=0.02, arrowsize=4pt]{->}(6.5,1.4){1.5}{210}{330}
\rput(6.5,-0.5){$A_{\Omega,0}$}

\rput(9.1,1.1){$\Omega_-$}
\rput(10.3,3.2){$\Omega_{+0}$}

\psline[linewidth=0.04](3,0)(1,2) 
\psline[linewidth=0.04](3,0)(3,5)
\psline[linewidth=0.04](3,0)(6,3)
\psline[linewidth=0.04](3,5)(6,3)
\psline[linewidth=0.04](3,5)(1,2)

\psline[linewidth=0.04](11.535,0.465)(8.707,0.465) 
\psline[linewidth=0.04](11.535,0.465)(8,4)
\psline[linewidth=0.04](11.535,0.465)(11.535,4.707)
\psline[linewidth=0.04](8,4)(11.535,4.707)
\psline[linewidth=0.04](8,4)(8.707,0.465)
\psdots[dotsize=4pt](8,4)(11.535,0.465)

\pspolygon[linewidth=0.03, linestyle=dashed](9.3035, 3.5535)(10.8535, 2.0535)(10.835, 3.835)

\psline[linewidth=0.02, arrowsize=4pt]{<->}(8.5,3.5)(8,4)\rput(8.5,3.9){$h^s$}

\rput(2,2){$h_\Omega\cdot K_-$}
\rput(4.3,2.8){$h_\Omega\cdot K_+$}
\rput(2.5,-0.5){$(\mathbf{0}, 0)$}
\rput(2.5, 5.5){$(\mathbf{0}, h_{\Omega})$}

\psdots[dotsize=4pt](3,0)(3,5)
\psbrace[linewidth=0.02, rot=135, nodesepA=-9pt](8.5,3.5)(11.035,0.965)
{$f_{\Omega 0}$}

\end{pspicture}
\caption{The transformation of the reference subdomain pair, the 
interior domain $\Omega_{+0}$ and the interior face $f_{\Omega 0}$ in the 
2-dimensional case.}
\label{K_figure}
\end{figure}
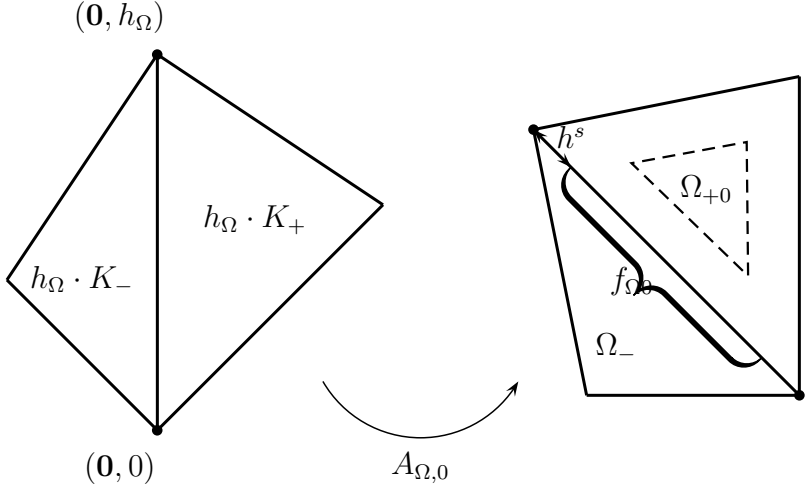

To extend the standard scaling arguments we first define a reference set 
$\mck$ of neighboring simplex pairs $(K_+, K_-)$ having the interelement 
face $f=\bar K_+\cap\bar K_-$ such that the following conditions hold:
\begin{itemize}
\item
$f\subset 0\times \er^{d-1}$ and one vertex of $f$ is 
$\mathbf{0}\in\er^d$
\item
the maximum edge-length of $f$ is \emph{one}
\item
$K_+$ and $K_-$ satisfy the condition on non-degeneracy.
\end{itemize}
Then for any neighboring subdomains $\Omega_+, \Omega_-\in \mct_h$ 
there is a pair $(K_+, K_-)\in\mck$ and an affine linear map 
$A_\Omega: K_+\cup K_-\to \Omega_+\cup \Omega_-$ with 
$A_\Omega(K_+) = \Omega_+$ and $A_\Omega(K_-) = \Omega_-$, moreover
\bee\label{basic_scale}
A_\Omega(\xx) = A_{\Omega,0}(h_\Omega \xx),
\eee 
where $h_\Omega$ denotes the maximum edge length of $f_\Omega$ 
and $A_{\Omega,0}$ is an isometry; see also Fig. \ref{K_figure}.
 
Accordingly, for any $v\in \mathbb{P}_{h,\kk}(\Omega_+\cup\Omega_-)$
the function $v_0:= v\circ A_{\Omega} \in \mathbb{P}_{h,\kk}(K_+\cup K_-)$,
moreover, using \eqref{basic_scale} the following equalities are valid: 
\bee\label{basic_scale2}
\jump{v_0}(\xx) = \jump{v} (A_{\Omega}(\xx))\quad\textrm{and}\quad
\eta_{h_0} * v_0 (\xx) = \eta_{h_0h^{\frac{1}{s}}} * v (A_\Omega \xx),
\eee  
whenever the operation $\eta_{h_0} *$ makes sense.

We also use the notation 
$h_\Omega\cdot K_\pm = \left\{h_\Omega\xx: \xx\in K_\pm \right\}$ and 
similarly $h_\Omega\cdot f$ and introduce the interior domain 
$\Omega_{j0} = \{\xx\in \Omega_j: B(\xx, h^s)\subset\Omega_j\}$ and
the interior face $f_0\subset f$ similarly.

The space $\textrm{BV}(\Omega)$ of real valued functions on $\Omega$
with bounded variations is defined with 
$$
\textrm{BV}(\Omega) = \left\{ 
u:\Omega\to\er: \sup_{\substack{\pphi\in [C_c^1(\Omega)]^d\\ 
                      \|\pphi\|_\infty = 1}} 
\int_\Omega u \nabla\cdot\pphi:= |u|_{\textrm{BV}}<\infty
\right\}
$$
and is equipped with the seminorm $|\cdot|_{\BV}$, where 
$\|\cdot\|_\infty$ denotes the maximum norm on $C_c^1(\Omega)$.
This seminorm can also be given as  
$$
|\cdot|_{\BV} = \int_{\Omega} \:\mathrm{d}|\pa u|,
$$
where $|\pa u|$ is the Radon measure generated by the distributional
derivative of $u$. 


The dual pairing between a distribution $S$ and a test function $\phi$ 
denoted using angle brackets: $\langle S, \phi\rangle$. 

In the estimates, the notation $g_1 \lm g_2$ means the existence of a 
constant $C$ - which does not depend on
the mesh parameter but possibly on the local polynomial degree - such 
that $g_1 \le C\cdot g_2$.  We also use the notation $g_1 \sim g_2$ 
provided that both $g_1 \lm g_2$ and $g_2 \lm g_1$ are satisfied. 

\section{Results}
The basic idea of the present analysis is to find a smoothed dG
approximation immediately. In this case, in the background we can 
compute with discontinuous basis functions in $\mathbb{P}_{h,\kk}$ and 
still have the freedom to choose them independently on the neighboring 
subdomains. 
On the other hand, as we compute conforming approximations, we can use 
the entire armory of the classical finite element analysis. 

The smoothed (or averaged) dG approximation consists of finding  
$\eta_h*u_h\in \mbp_{h,\kk,s}$
such that for all $\eta_h*v_h\in \mbp_{h,\kk,s}$
we have 
\bee\label{main}
a_\eta(u_h, v_h) := a_\eta^+(\eta_h*u_h, \eta_h*v_h):=
(\nabla (\eta_h*u_h), \nabla(\eta_h*v_h)) = (g_0, \eta_h*v_h), 
\eee
where the bilinear forms 
$a_\eta: \mbp_{h,\kk}\times \mbp_{h,\kk}\to\er$ and 
$a_\eta^+: \mbp_{h,\kk,s}\times \mbp_{h,\kk,s}\er$ are defined by 
$\eqref{main}$ and $g_0$ denotes the zero extension of $g$ to $\Omega_h$. 
Whenever the spaces $\mbp_{h,\kk,s}\not\subset H_0^1(\Omega)$ we call
the method $H^1$-conforming since each space is in $H_0^1(\Omega_h)$.  

We make use of the following inequalities, which can be proved using 
simple scaling arguments.
\begin{prop}
We have the following inequalities:
\bee\label{skala1}
\max_{B(\nnull, h^s)} |u| \sim h^{-\frac{sd}{2}} \|u\|_{B(\nnull, h^s)}
\qquad\forall u\in\mathbb{P}_{k}(B(\nnull,h^s)),
\eee
\bee\label{maxuv_and_l1norm}
\max_f\jump{u}\lm h^{1-d} \int_f |\jump{u}|
\qquad\forall \jump{u}\in\mathbb{P}_{k}(f),
\eee
\bee\label{max_nabla_sq_2_norm}
\max_K |\nabla^2 u|\lm
h^{-\frac{d}{2}-2} \|u\|_K
\qquad\forall u\in\mathbb{P}_{k}(K),
\eee
\bee\label{max_nabla2_l1norm}
\max_{f} \nabla_{f}^2 \jump{u} \lm h^{-d-1} \int_f |\jump{u}|  
\qquad\forall \jump{u}\in\mathbb{P}_{k}(f),
\eee
\bee\label{skala2}
\|\nabla u\|_{B(\nnull, h^s)} 
\lm h^{\frac{(s-1)d}{2}} \|\nabla u\|_{B(\nnull, h)}
\lm h^{-1} h^{\frac{(s-1)d}{2}} \|u\|_{B(\nnull, h)}\qquad\forall
u\in\mathbb{P}_k(B(\nnull,h)).\qquad\square
\eee
\end{prop}
We need also an estimate between the discontinuous function $\nabla_h u$ 
and its local average $\eta_h*\nabla_h u$ with a convergence rate depending 
on $h$. 
For this a Taylor expansion is developed about all $\xx\in\Omega_{j0}$ 
giving for an arbitrary $\yy\in\Omega_j$ that 
\bee\label{taylor_on_omegaj0}
u(\yy) = u(\xx) + \nabla u(\xx)\cdot(\yy-\xx) 
+ \frac{1}{2}\nabla^2 u(\xxi_\yy)(\yy-\xx)\cdot(\yy-\xx)   
\eee 
for some $\xxi_\yy$ in the section $(\xx, \yy)$. 
Integrating both sides over $B(\xx,h^s)$ yields
\bees
B_{h^s,d} \cdot (\eta_h* u(\xx)) = 
B_{h^s,d} \cdot u(\xx) + \int_{B(\xx,h^s)} 
\frac{1}{2}\nabla^2 u (\xxi_\yy)(\yy-\xx)\cdot(\yy-\xx)   
\eees 
and therefore
\bee\label{taylor_on_omegaj}
\eta_h* u(\xx) - u(\xx) =
\frac{1}{2\cdot B_{h^s,d}} 
\int_{B(\nnull,h^s)} \nabla^2 u(\xxi_\yy)|\yy|^2 \dyy.   
\eee 

\begin{prop}\label{prop_aver_approx}
For all $u\in \mathbb{P}_{h,\kk}$ and subdomain $\Omega_j$ we have 
\bees
\|\nabla_h u - \eta_h * \nabla_h u\|_{\Omega_j} 
\lm h^{\frac{s-1}{2}} \|\nabla_h u\|_{\tilde\Omega_j}.
\eees
\end{prop}
\emph{Proof:}
We first use the triangle inequality 
\bee\label{triang_for_conv}
\bali
&\|\nabla_h u - \eta_h * \nabla_h u\|_{\Omega_j} 
\le
\|\nabla_h u - \eta_h * \nabla_h u\|_{\Omega_{j0}} 
+
\|\nabla_h u - \eta_h * \nabla_h u\|_{\Omega_j\setminus\Omega_{j0}}\\
&
\le
  \|\nabla_h u - \eta_h * \nabla_h u\|_{\Omega_{j0}} 
+ \|\nabla_h u\|_{\Omega_j\setminus\Omega_{j0}} 
+ \|\eta_h * \nabla_h u\|_{\Omega_j\setminus\Omega_{j0}},
\eali
\eee 
where the contributions are estimated separately.
We obviously have the estimate 
$\lambda_d (\Omega_j\setminus\Omega_{j0})\lm h^s h_{\Omega_j}^{d-1}$ 
such that a simple scaling argument gives
\bee\label{first_part}
\|\nabla_h u\|_{\Omega_j\setminus\Omega_{j0}} \lm
h^{\frac{s-1}{2}} \|\nabla_h u\|_{\Omega_j}.
\eee
This also implies, using \eqref{skala2} in the second line with $s=1$   
that
\bee\label{last_part}
\bali
& \|\eta_h * \nabla_h u\|_{\Omega_j\setminus\Omega_{j0}}^2
\le
\lambda_d (\Omega_j\setminus\Omega_{j0})
\max_{\tilde\Omega_j} |\nabla_h u|^2\\
&\le
h^{d-1}_{\Omega_j}h^s\max_{\tilde\Omega_j} |\nabla_h u|^2
\le
h_{\Omega_j}^{d-1}h^s h_{\Omega_j}^{-d}\|\nabla_h u\|^2_{\tilde\Omega_j} 
= h^{s-1}\|\nabla_h u\|_{\tilde\Omega_j}^2.
\eali
\eee
Finally, combining the inequalities in \eqref{taylor_on_omegaj} and 
\eqref{max_nabla_sq_2_norm} we arrive at the estimate
$$
\bali
&|\nabla_h u - \eta_h * \nabla_h u|_{\Omega_{j0}}|
\le
\frac{1}{2\cdot\lambda_d(B(\xx, h^s))} 
\max_{\yy\in\Omega_j}|\nabla^3 u(\yy)| \int_{B(\nnull,h^s)} |\yy|^2 \dyy\\
&\lm
h^{-sd} \max_{\yy\in\Omega_j}|\nabla^3 u(\yy)| h^{s(d+2)}
\lm
h^{2s} h^{-\frac{d}{2}-2}\|\nabla u\|_{\Omega_{j}}.   
\eali
$$
Therefore, using \eqref{skala1} we obtain
\bee\label{third_part}
\|\nabla_h u - \eta_h * \nabla_h u\|_{\Omega_{j0}}
\le
h^{2s-2}\|\nabla_h u\|_{\Omega_{j}}.
\eee
The estimates \eqref{first_part}, \eqref{last_part} and
\eqref{third_part} with \eqref{triang_for_conv} imply then
the inequality in the proposition.
\eproof

\emph{Remark:}
For functions $v\in C^2(\er^d)$ one can easily estimate the difference in 
Proposition \ref{prop_aver_approx}. Moreover, it turns out that the convergence rate of the difference $\int_{\er^d}|\eta_h*v|^2 - |v|^2$ characterizes the Sobolev space $H^1(\er^n)$, see \cite{peletier07}.

The chief problem in the estimations with convolution terms is that 
the scaling arguments can not be applied in a straightforward way.
Whenever we use polynomial spaces the function space 
$\{\eta_h * v: v\in\mathbb{P}_{h,\kk}, 0<h<h_0\}$
is infinite dimensional, which makes the following proofs non-trivial.

\begin{prop}\label{prop2}
There exists $h_0>0$ such that for all $h$ with $h^{1-\frac{1}{s}}<h_0$ 
and $v\in \mathbb{P}_{h,\kk}(\Omega_+\cup \Omega_-)$ we have
\bee\label{prop2_est}
\int_{f_\Omega} |\jump{v}| \lm 
\int_{\Omega_+\cup \Omega_-} |\nabla (\eta_h* v)|
\eee 
and for $s\ge \frac{3}{2}$
\bee\label{prop3_est}
\int_{f_\Omega} |\jump{v}| \lm
h^{\frac{d}{2}} 
\sqrt{\int_{\Omega_+\cup \Omega_-} |\nabla (\eta_h* v)|^2}.
\eee 
\end{prop}
The corresponding proof is postponed to the Appendix.

\subsection{The bilinear form}
To give the bilinear form \eqref{main} in a more explicit form, we first 
need some identities for distributional derivatives.

We first decompose the gradient 
of a function $u\in \mathbb{P}_{h,\kk}(\Omega)$ as follows.
\begin{lem}\label{lem1}
For all $u\in \mathbb{P}_{h,\kk}(\Omega)$ we have 
$$
\nabla u = \nabla_h u + \jump{u}_{\mcd}
$$ 
in the sense of distributions, i.e. 
$\jump{u}_\mcd\in [\mcd^3(\Omega)]^*$ is a distribution with
\[
\langle\jump{u}_\mcd, \pphi\rangle 
= -\sum_{f\in\mcf} \int_f \jump{u}_f \cdot \pphi
:= -\int_{\mcf} \jump{u} \cdot \pphi = -(\jump{u}, \pphi)_\mcf.
\] 
\end{lem}
\proof
Obviously, for all $\pphi\in[\mcd(\Omega)]^3$ we have 
\[
\bali
&
\langle \nabla u, \pphi\rangle = 
- \langle u, \nabla\cdot \pphi\rangle =
- \sum_{\Omega_j\in\mathcal{T}_h}\int_{\Omega_j} u \nabla\cdot \pphi = 
\sum_{\Omega_j\in\mathcal{T}_h}\int_{\Omega_j} \nabla u \cdot\pphi - 
\sum_{\Omega_j\in\mathcal{T}_h}\int_{\pa\Omega_j} u|_{\Omega_j} \nnu_j\cdot\pphi\\
&=
 \sum_{\Omega_j\in\mathcal{T}_h}\int_{\Omega_j} \nabla u \cdot\pphi - 
\sum_{f\in\mcf}\int_{f} \jump{u}_f \cdot \pphi|_f,
\eali
\]
which proves the statement.
\eproof

\emph{Remarks:} The decomposition in Lemma \ref{lem1} is indeed
a Lebesgue decomposition \cite{halmos50} of the Radon measure 
corresponding to the distributional derivative $\nabla u$, which  
can be considered as a special case of the one in \cite{szucs13}. The role 
of the jump terms in this context in analyzed in \cite{attouch06}, Section 10.\\
The symbol $\jump{\cdot}_\mcd$ can be understood both as a distribution 
supported on the interelement faces and the singular measure in the 
corresponding Lebesgue decomposition. The connection between $\jump{u}_\mcd$
with classical function $\jump{u}$ is highlighted in Lemma \ref{lem1}.\\
The negative sign is a weakness of the conventional notation. 
This is already transparent in the one-dimensional case: whenever 
the Heaviside step function $H:\er\to\er$ is increasing, by definition 
we have $\jump{H}(0)=-1$.


For the consecutive derivations we need also an identity regarding the 
convolution of distributions. 

\begin{lem}\label{lem2}
For all $u\in \mathbb{P}_{h, \kk}$ the convolution $\eta_h*\jump{u}_{\mcd}$
is regular, which will be identified with the corresponding
locally integrable function. With this, for all 
bounded function $\ww: \Omega\to\er^3$ we have  
\[
\langle \eta_h*\jump{u}_\mcd, \ww\rangle = 
(\jump{u}, \eta_h*\ww)_\mcf.
\]
\end{lem}
\proof
Since both $\eta_h$ and $\jump{u}$ are compactly supported, 
we get by definition (see \cite{hirsch99}, Definition 2.1) and 
by Lemma \ref{lem1} that for each $\pphi\in [C_0^\infty(\Omega)]^3$ 
the following equality is valid:
\bee\label{derive_lem2}
\bali
&\langle \eta_h*\jump{u}_\mcd, \pphi\rangle 
= \langle\jump{u}_\mcd, \yy\to \eta_h(\xx\to\pphi(\xx+\yy))\rangle 
= \langle\jump{u}_\mcd, \yy\to \int_{\er^d} \eta_h(\xx)\pphi(\xx+\yy)\dxx\rangle\\
&=-\int_\mcf \jump{u}(\yy) \int_{\er^d} \eta_h(\xx)\pphi(\xx+\yy)\dxx\dyy
=-\int_\mcf \jump{u}(\yy) \int_{\er^d} \eta_h(\zz-\yy)\pphi(\zz)\dzz\dyy\\
&=-\int_\mcf \jump{u}(\yy) \int_{\er^d} \eta_h(\yy-\zz)\pphi(\zz)\dzz\dyy
=-\int_\mcf \jump{u}(\yy)\;\eta_h*\pphi(\yy) \dyy
=-(\jump{u}, \eta_h*\pphi)_\mcf.
\eali
\eee
On the other hand, according to \cite{hirsch99}, page 337, 
Exercise 10, $\eta_h*\jump{u}$ is locally integrable
such that the statement of the lemma is valid for all
bounded functions $\ww$ as it was stated.\quad $\square$

Then we get as an obvious consequence the following.
\begin{cor}
The bilinear form $a_\eta$ can be rewritten as
\bee\label{newform}
\bali
&a_\eta(u,v) \\
&= (\eta_h*\nabla_h u, \eta_h*\nabla_h v) + 
\langle\eta_h*\nabla_h u, \eta_h*\jump{v}_\mcd\rangle +
\langle\eta_h*\nabla_h v, \eta_h*\jump{u}_\mcd\rangle +
(\eta_h*\jump{u}, \eta_h*\jump{v})\\
&= (\eta_h*\nabla_h u, \eta_h*\nabla_h v) - 
(\eta_h*\eta_h*\nabla_h u, \jump{v})_\mcf -
(\eta_h*\eta_h*\nabla_h v, \jump{u})_\mcf +
(\eta_h*\jump{u}, \eta_h*\jump{v}).
\eali
\eee
\end{cor}
Note that the first line is related to the lifted forms 
of the dG methods as each scalar product corresponds to a volume 
integral. On the other hand, the second and third terms in the 
second line are integrals which can be computed on faces 
according to the second line in \eqref{derive_lem2}. 


\section{Comparison with the IP bilinear form}\label{sect_4}

We compare our bilinear form \eqref{newform} with the IP bilinear form 
\eqref{IP_form} componentwise.  

The first lemma quantifies the difference of the first terms.
\begin{lem}\label{diff_first_terms_lem}
For all $u,v \in\mathbb{P}_{h,\kk}(\Omega)$ we have
\bees
|(\nabla_h u, \nabla_h v) - (\eta_h * \nabla_h u, \eta_h * \nabla_h v)|
\le
h^{s-1} \|\eta_h * \nabla_h u\|  \|\eta_h * \nabla_h v\|.
\eees
\end{lem}
\emph{Proof:}
We obviously have
\bee\label{diff_first_terms_1}
\bali
&|(\nabla_h u, \nabla_h v) - (\eta_h * \nabla_h u, \eta_h * \nabla_h v)|\\
&\le
|(\nabla_h u - \eta_h * \nabla_h u, \nabla_h v) + 
 (\eta_h * \nabla_h u, \nabla_h v - \eta_h * \nabla_h v) |\\
&\le
\|\nabla_h u - \eta_h * \nabla_h u\| \|\nabla_h v\| + 
 \|\eta_h * \nabla_h u\| \|\nabla_h v - \eta_h * \nabla_h v\|.
\eali
\eee
Also, application of the estimate in Proposition \ref{prop_aver_approx} 
and a simple scaling argument implies for each subdomain $\Omega_j$ that
\bees
\bali
&\|\nabla_h u|_{\Omega_{j0}}\| \le 
\|\nabla_h u - \eta_h*\nabla_h u\|_{\Omega_{j0}} + 
\|\eta_h*\nabla_h u\|_{\Omega_{j0}}\\
&\le
h^{\frac{s-1}{2}} \|\nabla_h u\|_{\Omega_{j}} + \|\eta_h*\nabla_h u\|_{\Omega_{j0}}
\lm
h^{\frac{s-1}{2}} \|\nabla_h u\|_{\Omega_{j0}} + \|\eta_h*\nabla_h u\|_{\Omega_{j0}}
\eali
\eees
and therefore,
\bees
\|\nabla_h u\|_{\Omega_{j0}} \lm \|\eta_h*\nabla_h u\|_{\Omega_{j0}},
\eees
which can be used to obtain the following inequality:
\bee\label{normal_vs_conv_norm}
\|\nabla_h u\|_{\Omega_j} \lm \|\nabla_h u\|_{\Omega_{j0}} \lm 
 \|\eta_h*\nabla_h u\|_{\Omega_{j0}}\le \|\eta_h*\nabla u\|_{\Omega_{j}}.  
\eee
Therefore, using again Proposition \ref{prop_aver_approx} we also have
\bee\label{normal_vs_conv_norm2}
\|\nabla_h u - \eta_h * \nabla_h u\|_{\Omega_j} 
\lm h^{\frac{s-1}{2}} \|\eta_h * \nabla u\|_{\Omega_j}. 
\eee
Taking the square of \eqref{normal_vs_conv_norm} and 
\eqref{normal_vs_conv_norm2} for each index $j$ and summing 
them we have 
$$
 \|\nabla_h u\|\lm \|\eta_h*\nabla u\|\quad\textrm{and}\quad
\|\nabla_h u - \eta_h * \nabla_h u\|
\lm h^{\frac{s-1}{2}} \|\eta_h * \nabla u\| 
$$
which can be used in \eqref{diff_first_terms_1} to obtain 
$$
\bali
&|(\nabla_h u, \nabla_h v) - (\eta_h * \nabla_h u, \eta_h * \nabla_h v)|\\
&\lm
h^{\frac{s-1}{2}} \|\eta_h * \nabla u\| \|\eta_h*\nabla v\| +
h^{\frac{s-1}{2}} \|\eta_h * \nabla v\| \|\eta_h*\nabla u\| 
\eali
$$
as stated in the lemma.
\eproof

To compare the second and third terms in \eqref{newform} and \eqref{IP_form}
we use the notation in Fig. \ref{K_figure} and the corresponding explanation.

To analyze the average of the approximations
we use the following statement on integral means. 
\begin{prop}\label{prop_seged}
For each $u\in \mathbb{P}_{h,\kk}(\Omega_-\cup\Omega_+)$ and 
$\xx\in f_\Omega$ with $B(\xx,2h^s)\subset \Omega_-\cup\Omega_+$ there exist 
$\bar\xx_- \in B_-(\xx, 2h^s)$ and $\bar\xx_+ \in B_+(\xx, 2h^s)$ such that
$$
u(\bar\xx_-) = 
2 \int_{B_-(\nnull, 2h^s)} u(\xx-\zz) \cdot \eta_h*\eta_h (\zz) \dzz 
$$
and similarly,
$$
u(\bar\xx_+) = 
2 \int_{B_+(\nnull, 2h^s)} u(\xx-\zz) \cdot \eta_h*\eta_h (\zz) \dzz, 
$$
where $B_-(\nnull, 2h^s)$ and $B_+(\nnull, 2h^s)$ denote the half-ball with 
non-positive and non-negative first coordinates, respectively.  
\end{prop}
The proof is postponed to the Appendix.

\begin{prop}\label{prop1}
For all $u\in\mathbb{P}_{h,\kk}(\Omega_-\cup\Omega_+)$ we have the
following inequality:
\bee\label{lem_11_ineq}
\max_{
\substack{\xx\in f_\Omega
        \\B(\xx, 2h^s)\subset\Omega_+\cup\bar\Omega_-}}
|\eta_h*\eta_h*\nabla_h u (\xx) - \aver{\nabla_h u}(\xx)| \lm 
h^{s-\frac{d}{2}-1}
\|\nabla_h u\|_{B(\Omega_+\cup\Omega_-)}. 
\eee
\end{prop}
\proof
Using the result of Proposition \ref{prop_seged} we rewrite the difference 
on the left hand side of \eqref{lem_11_ineq} as follows:
\bee\label{deriv1_34terms}
\bali
&
\eta_h*\eta_h*\nabla_h u (\xx) - \aver{\nabla_h u}(\xx)\\
& =
\frac{1}{2} 
\left(
2 \int_{B_-(\nnull, 2h^s)}\nabla_h u(\xx-\zz)\cdot \eta_h*\eta_h (\zz)\dzz 
+
2 \int_{B_+(\nnull, 2h^s)}\nabla_h u(\xx-\zz)\cdot \eta_h*\eta_h (\zz)\dzz\right.\\
\qquad&\qquad\left.
- \nabla_h u(\xx_-) -\nabla_h u(\xx_+)
\right) \\
& =
\frac{1}{2} 
\left(
\nabla_h u(\bar\xx_-) -\nabla_h u(\xx_-)
+ \nabla_h u(\bar\xx_+) -\nabla_h u(\xx_+)
\right).
\eali
\eee
We use then the estimate 
\bees
|\nabla_h u(\bar\xx_-) - \nabla_h u(\xx_-)| 
\le 
2h^s\cdot \sup_{\zz\in (\xx_-, (0-,\yy))}  \|\nabla^2_h u(\zz)\|
\eees
in \eqref{deriv1_34terms} to see that
\bees
\bali
&\max_{\substack{\xx\in f_\Omega
        \\B(\xx, 2h^s)\subset\Omega_+\cup\bar\Omega_-}}
|\eta_h*\eta_h*\nabla_h u (\xx) - \aver{\nabla_h u}(\xx)|\\
& \le
\frac{1}{2}\cdot 2h^s\cdot 
\left(
\max_{\zz\in f_\Omega\dotplus B(\nnull, 2h^s)}  \|\nabla^2_h u(\zz)\| 
+ 
\max_{\zz\in f_\Omega\dotplus B(\nnull, 2h^s)}  \|\nabla^2_h u(\zz)\|
\right)\\
&=
2h^s
\max_{\zz\in f_\Omega\dotplus B(\nnull, 2h^s)}  \|\nabla^2_h u(\zz)\|.
\eali
\eees
The last term here can be estimated using scaling arguments as
\bees
\bali
&
2h^s\max_{\zz\in f_\Omega\dotplus B(\nnull, 2h^s)}  
\|\nabla^2_h u(\zz)\|
\lm
h^s \sqrt{h^{-s}h_\Omega^{1-d}}\|\nabla^2_h u\|_{f_\Omega\dotplus B(\nnull, 2h^s)}\\
&
\lm
h^{\frac{1+s-d}{2}} \|\nabla^2_h u\|_{f_\Omega\dotplus B(\nnull, 2h^s)}
\lm
h^{\frac{1+s-d}{2}} h^{\frac{s-1}{2}} \|\nabla^2_h u\|_{\Omega_+\cup\Omega_-}
\lm
h^{s-\frac{d}{2}-1} \|\nabla_h u\|_{\Omega_+\cup\Omega_-},
\eali
\eees
which proves the statement of the proposition. \quad $\square$

We can now relate the third and second terms in the proposed  
bilinear form \eqref{newform} and the IP bilinear form.
\begin{lem}\label{23terms}
For arbitrary $u,v\in \mathbb{P}_{h,\kk}$ we have the following 
inequality:
\bees
|(\eta_h*\eta_h*\nabla_h u, \jump{v})_\mcf
- (\aver{\nabla_h u}, \jump{v})_\mcf|
\lm h^{s-1}  \|\nabla (\eta_h* u)\| \|\nabla(\eta_h* v)\|.
\eees
\end{lem}
\emph{Proof:}
Using the result of Proposition \ref{prop1} and Proposition \ref{prop2}
we obtain the following estimation on the interelement face $f_\Omega$
between $\Omega_+$ and $\Omega_-$: 
\bees
\bali
&|(\eta_h*\eta_h*\nabla_h u, \jump{v})_{f_\Omega}
- (\aver{\nabla_h u}, \jump{v})_{f_\Omega}|=
|(\eta_h*\eta_h*\nabla_h u - \aver{\nabla_h u}, \jump{v})_{f_\Omega}|\\
&
\le
\max_{\xx\in f_\Omega}|(\eta_h*\eta_h*\nabla_h u - \aver{\nabla_h u}(\xx)|
\int_{f_\Omega}|\jump{v}|\\
&
\lm
\max_{\substack{\xx\in f_\Omega
      \\B(\xx, 2h^s)\subset\Omega_+\cup\bar\Omega_-}}|
(\eta_h*\eta_h*\nabla_h u - \aver{\nabla_h u}(\xx)|
\int_{f_\Omega}|\jump{v}|\\
&
\le
h^{s-\frac{d}{2}-1}\left(\frac{h}{h_\Omega}\right)^{\frac{d}{2}}  
\|\nabla_h u\|_{\Omega_+\cup \Omega_-} \int_f |\jump{v}|
\le 
h^{s-\frac{d}{2}-1}\left(\frac{h}{h_\Omega}\right)^{\frac{d}{2}}
 \|\nabla_h u\|_{\Omega_+\cup \Omega_-}
\int_{\Omega_+\cup \Omega_-} |\nabla (\eta_h* v)|\\
&\lm
h^{s-\frac{d}{2}-1}\left(\frac{h}{h_\Omega}\right)^{\frac{d}{2}}
\|\nabla_h (\eta_h* u)\|_{\Omega_+\cup\Omega_-}
h_\Omega^{\frac{d}{2}}
\|\nabla (\eta_h* v)\|_{\Omega_+\cup\Omega_-}\\
&\lm
h^{s-1}
\|\nabla_h (\eta_h* u)\|_{\Omega_+\cup\Omega_-}
\|\nabla (\eta_h* v)\|_{\Omega_+\cup\Omega_-}.
\eali
\eees
Summing up these inequalities for each interelement face $f_\Omega$ 
and using the discrete Cauchy--Schwarz inequality result in the estimate
$$
\bali
&|\sum_{f_\Omega\in\mcf}
(\eta_h*\eta_h*\nabla_h u, \jump{v})_{f_\Omega}
- \sum_{f_\Omega\in\mcf} (\aver{\nabla_h u}, \jump{v})_{f_\Omega}|
\le
\sum_{f_\Omega\in\mcf}
|(\eta_h*\eta_h*\nabla_h u, \jump{v})_{f_\Omega}
- (\aver{\nabla_h u}, \jump{v})_{f_\Omega}|\\
&\le
\sum_{f_\Omega\in\mcf}
h^{s-1}  
\|\nabla_h (\eta_h* u)\|_{\Omega_+\cup\Omega_-}
\|\nabla_h (\eta_h* v)\|_{\Omega_+\cup\Omega_-}\\
&\lm
h^{s-1}
\sqrt {\sum_{f_\Omega\in\mcf} \|\nabla (\eta_h* u)\|^2_{\Omega_+\cup\Omega_-}}
\sqrt {\sum_{f_\Omega\in\mcf} \|\nabla (\eta_h* v)\|^2_{\Omega_+\cup\Omega_-}}
\lm
h^{s-1} \|\nabla (\eta_h* u)\|^2  \|\nabla (\eta_h* v)\|^2
\eali
$$
as stated in the lemma.\quad$\square$

To relate the last term in \eqref{newform} with the penalty term in
the IP bilinear form, we rewrite the locally integrable function
$\eta_h*\jump{v}$ (see Lemma \ref{lem2}) in a more explicit form. 
\begin{lem}\label{compute_conv}
For each $v\in\mathbb{P}_{\hh,k}$ and $f\in\mcf$ the following  
identity is valid:
\bee\label{expand_conv}
\eta_h*\jump{v}_f (\xx) 
=
\int_f \eta_h (\xx-\yy) \jump{v}_f(\yy)\:\mathrm{d}\yy. 
\eee
\end{lem}
This result can also serve as a good argument why did we apply 
the same notation for the convolution corresponding to the jump 
of $v$ and the jump function. Since the proof is a bit technical 
it is postponed to the appendix.

\begin{figure}
\begin{pspicture}(0,-1)(10, 3)
\psframe[linewidth=0cm,fillcolor=lightgray, fillstyle=solid](1.99,0)(9.02,2)
\psarc[linewidth=0cm,fillcolor=lightgray, fillstyle=solid](2,1){1}{90}{270}
\psarc[linewidth=0cm,fillcolor=lightgray, fillstyle=solid](9,1){1}{270}{90}

\psline[linewidth=0.04](2,1)(9,1)
\psline[linewidth=0.03, linestyle=dotted](1.2,1.6)(9.8,1.6)
\psline[linewidth=0.03, linestyle=dotted](1.2,0.4)(9.8,0.4)
\psline[linewidth=0.03, linestyle=dashed](3,0.4)(8,0.4)
\psline[linewidth=0.02](3.5,1)(4.3,1.6)
\psline[linewidth=0.02](5.1,1)(4.3,1.6)

\psline[linewidth=0.02]{|<->|}(9,1)(10,1)\rput(9.5,0.8){$h^s$}
\psline[linewidth=0.02]{|<->|}(2,1.6)(2,1)\rput(2.2,1.3){$r$}
\psline[linewidth=0.02]{|<->|}(2,0.4)(2,1)\rput(2.2,0.7){$r$}

\rput(4.3,1.8){$\mathbf{x}$}
\psdots[dotsize=4pt](4.3,1.6)
\psbrace[linewidth=0.02, ref=lC,rot=90, nodesepB=-8pt, nodesepA=-4pt](3.5,1)(5.1,1){$f_{\mathbf{x},r}$}

\psframe[linewidth=0.02](11.4,0)(13.8,2)

\psline[linewidth=0.03, linestyle=dotted](11.5,1.6)(12.5,1.6) 
\rput(13.2,1.6){$f\otimes r$}
\psline[linewidth=0.04, linestyle=solid](11.5,1)(12.5,1) 
\rput(13.2,1){$f$}
\psline[linewidth=0.03, linestyle=dashed](11.5,0.4)(12.5,0.4) 
\rput(13.2,0.4){$f_0\dot{+} r\textrm{\boldmath$\nu$}$}
\end{pspicture}
\caption{Interelement face $f$ with the support of 
$\eta_h * \jump{u}_f$ (shaded) and the sets $f_{\xx,r}$, $f\otimes r$
and $f_0\dot{+}r\nnu$.}
\label{pict_ftrace}
\end{figure}
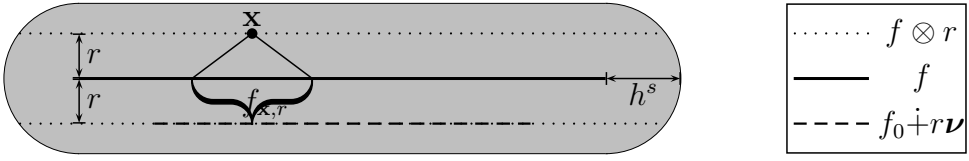

To analyze the right hand side of \eqref{expand_conv}, we introduce 
the following sets which are depicted in Figure \ref{pict_ftrace}.
\bees
f\otimes r = \{\xx\in \langle f\rangle \dotplus r \nnu_1 \cup
\langle f\rangle \dotplus r \nnu_2: d(\xx,f)\le h^s\}
\eees
and
\bees
f_0\otimes r =  (f_0 \dotplus r \nnu_1) \cup (f_0 \dotplus r \nnu_2).
\eees,
where $\langle f\rangle$ denotes the affine subspace generated by $f$.

Observe that $\eta_h*\jump{u}_f (\xx)$ can be nonzero if $\xx\in f\otimes r$ 
for some $r< h^s$ and then we use the notation $f_{\xx,r} = B(\xx,h^s)\cap f$,
which is a ball  in $\langle f \rangle$ centered at the projection of $\xx$
on $f$ with the radius $\sqrt{h^{2s}-r^2}$ such that 
$\lambda(f_{\xx,r}) = B_{\sqrt{h^{2s}-r^2}, d-1}$.

With these, we can rewrite \eqref{expand_conv} as 
\bees
\eta_h*\jump{u}_f (\xx) = 
\frac{1}{B_{h^s, d}}
\int_{f_{\xx,r}} \jump{u}. 
\eees  
In this way, using Lemma \ref{compute_conv} the integral in the last 
term of \eqref{newform} \emph{on a face $f$} can be rewritten as
\bee\label{exact_last}
(\eta_h*\jump{u}_f, \eta_h*\jump{v}_f)
= 
\frac{1}{[B_{h^s, d}]^2}
\int_{-h^s}^{h^s}
\int_{f\otimes r}
\int_{f_{\xx,r}} \jump{u} (\sss)\dss\int_{f_{\xx,r}} \jump{v} (\sss)
\dss\dxx\dr.
\eee
We intend to relate this term with the following:
\bee\label{appr_last}
\frac{1}{[B_{h^s,d}]^2} \int_{-h^s}^{h^s} \int_{f}
[B_{\sqrt{h^{2s}-r^2},d-1}]^2 \jump{u}\jump{v} (\xx) 
\dxx\dr.
\eee
To work with smooth functions, both in \eqref{exact_last} and 
\eqref{appr_last} we have to restrict the integrals on 
$f_0\otimes r$ and to $f_0$, respectively. 
Since $\lambda(f)\sim h^{d-1}$ and  
$\lambda(f\setminus f_0)\sim h^{d-2}h^s$, a scaling argument 
implies the following estimates:
\bee\label{fortriag1}
\bali
& I_1(r):=\\
&\left|\int_{f\otimes r}
\int_{f_{\xx,r}} \jump{u} (\sss)\dss\int_{f_{\xx,r}} \jump{v} (\sss)
\dss\dxx
-
\int_{f_0}
\int_{B(\xx,\sqrt{h^{2s}-r^2})} \jump{u} (s)\dss
\int_{B(\xx,\sqrt{h^{2s}-r^2})} \jump{v} (s)\dss \dxx\right| \\
&=
\left|\int_{f\otimes r}
\int_{f_{\xx,r}} \jump{u} (\sss)\dss\int_{f_{\xx,r}} \jump{v} (\sss)
\dss\dxx
-
\int_{f_0\otimes r}
\int_{f_{\xx,r}} \jump{u} (\sss)\dss\int_{f_{\xx,r}} \jump{v} (\sss)
\dss\dxx\right| \\
&\lm
\frac{h^{d-2}h^s}{h^{d-1}}
\int_{f\otimes r}
\left|\int_{f_{\xx,r}} \jump{u} (\sss)\dss\right|
\left|\int_{f_{\xx,r}} \jump{v} (\sss)\dss\right|
\dxx =
h^{s-1}
\int_{f\otimes r}
\left|\int_{f_{\xx,r}} \jump{u} (\sss)\dss\right|
\left|\int_{f_{\xx,r}} \jump{v} (\sss)\dss\right|
\dxx
\eali
\eee
and 
\bee\label{fortriag2}
\bali
& I_2(r):=
\left| 
\int_{f_0}
[B_{\sqrt{h^{2s}-r^2}, d-1}]^2 \jump{u} \jump{v} (\xx)\dxx
-
\int_{f}
[B_{\sqrt{h^{2s}-r^2}, d-1}]^2 \jump{u} \jump{v} (\xx)\dxx
\right|\\
&\lm
\frac{h^{d-2}h^s}{h^{d-1}}
\int_{f}
[B_{\sqrt{h^{2s}-r^2}, d-1}]^2 |\jump{u} \jump{v}| (\xx)\dxx
\lm
h^{s-1}(h^{2s}-r^2)^{d-1} 
\int_{f} |\jump{u} \jump{v}| (\xx)\dxx.
\eali
\eee
\emph{Remark:}
The estimation of $I_1$ in still valid if we use $f_{00}\subset f$  
with $\lambda(f\setminus f_{00})\sim h^{d-2}h^s$.

For the forthcoming computations, we also give the magnitude of
the following integrals:
\bee\label{magn1_int}
\int_{-h^s}^{h^s}(h^{2s}-r^2)^{d-1}\dr = \ordo(h^{2sd-s})
\eee
\bee\label{magn2_int}
\int_{B(\xx, \sqrt{h^{2s}-r^2})} |\sss-\xx|^2\dss 
= \ordo(\sqrt{h^{2s}-r^2}),
\eee
which can be verified with a straightforward computation.
\begin{lem}\label{est_last_aver_lem}
For all $u,v\in\mathbb{P}_{h,\kk}$ and $\Omega_+, \Omega_-\in\mathcal{T}_h$ 
we have the following inequality
\bee\label{est_last_aver}
\bali
&\left|(\eta_h*\jump{u}_f, \eta_h*\jump{v}_f)
-
\frac{1}{[B_{h^s,d}]^2} \int_{-h^s}^{h^s} \int_{f}
[B_{\sqrt{h^{2s}-r^2},d-1}]^2 \jump{u}\jump{v} (\xx) 
\dxx\dr\right|\\
&\le
h^{s-1}(1+h^{3s-d-2}) \|\nabla(\eta_{h}*u)\|_{\Omega_+\cup\Omega_-} 
\|\nabla(\eta_{h}*v)\|_{\Omega_+\cup\Omega_-}. 
\eali
\eee
\end{lem}
\emph{Proof:}
Using \eqref{exact_last} and a triangle inequality with 
\eqref{fortriag1} and \eqref{fortriag2} we have
\bee\label{4term_longest}
\bali
&\left|
(\eta_h*\jump{u}_f, \eta_h*\jump{v}_f) -
\frac{1}{[B_{h^s,d}]^2} \int_{-h^s}^{h^s} \int_{f}
[B_{\sqrt{h^{2s}-r^2},d-1}]^2 \jump{u}\jump{v} (\xx) 
\dxx\dr
\right|\\
&\le
\frac{1}{[B_{h^s,d}]^2} \int_{-h^s}^{h^s} I_1(r)\dr
+
\frac{1}{[B_{h^s,d}]^2} \int_{-h^s}^{h^s} I_2(r)\dr\\
&+
\frac{1}{[B_{h^s,d}]^2} 
\int_{-h^s}^{h^s}
\left|
\int_{f_0}
\int_{B(\xx,\sqrt{h^{2s}-r^2})} \jump{u} (\sss)\dss
\int_{B(\xx,\sqrt{h^{2s}-r^2})} \jump{v} (\sss)\dss\dxx\right.\\
\qquad\qquad &\qquad \left. - 
 \int_{f_0}
[B_{\sqrt{h^{2s}-r^2},d-1}]^2 \jump{u} \jump{v} (\xx)\dxx
\right|\dr\\
&\le
\frac{h^{s-1}}{[B_{h^s,d}]^2}
\int_{-h^s}^{h^s}
\int_{f\otimes r}
\left|\int_{f_{\xx,r}} \jump{u} (\sss)\dss\right|
\left|\int_{f_{\xx,r}} \jump{v} (\sss)\dss\right|\dxx\dr\\
&+
\frac{h^{s-1}}{[B_{h^s,d}]^2}
\int_{-h^s}^{h^s}
 \int_{f}
[B_{\sqrt{h^{2s}-r^2},d-1}]^2 |\jump{u} \jump{v}| (\xx)\dxx\dr\\
&+
\frac{1}{[B_{h^s,d}]^2} 
\int_{-h^s}^{h^s}
\left|
\int_{f_0}
\int_{B(\xx,\sqrt{h^{2s}-r^2})} \jump{u} (\sss)\dss
\int_{B(\xx,\sqrt{h^{2s}-r^2})} \jump{v} (\sss)\dss\right.\\
\qquad
&\left. - 
[B_{\sqrt{h^{2s}-r^2},d-1}]^2 \jump{u} \jump{v} (\xx)\dxx
\right|\dr.
\eali
\eee
The error terms here are estimated separately.\\
We first use \eqref{maxuv_and_l1norm} and \eqref{magn1_int} to obtain
\bee\label{1term_last}
\bali
&\frac{h^{s-1}}{[B_{h^s,d}]^2}\int_{-h^s}^{h^s}
\int_{f\otimes r}
\left|\int_{f_{\xx,r}} \jump{u} (\sss)\dss\right|
\left|\int_{f_{\xx,r}} \jump{v} (\sss)\dss\right|\dxx\dr\\
&\lm
h^{s-2sd-1} \int_{-h^s}^{h^s}
\lambda_{d-1}(f\otimes r) 
[B_{\sqrt{h^{2s}-r^2},d-1}]^2 \max_f\jump{u} \max_f\jump{v}\dr\\
&\lm
h^{s-2sd-1} 
\int_{-h^s}^{h^s}
h_\Omega^{d-1} (h^{2s}-r^2)^{d-1}\cdot 
h_\Omega^{1-d} \int_f |\jump{u}| 
h_\Omega^{1-d} \int_f |\jump{v}|\dr \\
&\le
h^{s-2sd-d}
\int_{-h^s}^{h^s} (h^{2s}-r^2)^{d-1}\dr 
 \int_f |\jump{u}| \int_f |\jump{v}|\\
&=
h^{s-2sd-d}h^{2sd-s}
\int_f |\jump{u}| \int_f |\jump{v}|
=
h^{-d}\int_f |\jump{u}| \int_f |\jump{v}|.
\eali
\eee
We proceed similarly for the second term in \eqref{4term_longest}:
\bees
\bali
&\frac{h^{s-1}}{[B_{h^s,d}]^2}\int_{-h^s}^{h^s} \int_{f}
[B_{\sqrt{h^{2s}-r^2},d-1}]^2 
\left|\jump{u} (\xx)\jump{v} (\xx)\right|\dxx\dr\\
&\lm
h^{s-2sd-1} h^{2sd-s}\int_{f}|\jump{u}\jump{v}|
\lm
h^{-1} h_\Omega^{1-d}
\int_{f}|\jump{u}| \int_{f}|\jump{v}|
\le
h^{-d}
\int_{f}|\jump{u}| \int_{f}|\jump{v}|.
\eali
\eees
We finally estimate the third term in \eqref{4term_longest}.
Using the expansion in \eqref{taylor_on_omegaj0} on $f_0$ with the 
surface gradient $\nabla _f \jump{u} := \nabla \jump{u}$
and integrating both sides on the ball $B(\xx, \sqrt{h^{2s}-r^2})$ implies
\bee\label{midp_aver}
\int_{B(\xx, \sqrt{h^{2s}-r^2})}\jump{u} (\sss)\dss = 
B_{\sqrt{h^{2s}-r^2},d-1} \jump{u}(\xx)
+ \frac{1}{2} \int_{B(\xx, \sqrt{h^{2s}-r^2})}
\nabla^2 \jump{u}(\xxi_\sss)(\sss-\xx)\cdot(\sss-\xx)\dss.  
\eee
Taking the product of \eqref{midp_aver} for $\jump{u}$ and $\jump{u}$ 
and using \eqref{magn2_int} and \eqref{max_nabla2_l1norm} we obtain
\bees
\bali
&
\left|
\int_{f_0}
\int_{B(\xx,\sqrt{h^{2s}-r^2})} \jump{u} (\sss)\dss
\int_{B(\xx,\sqrt{h^{2s}-r^2})} \jump{v} (\sss)\dss\dxx  - 
[B_{\sqrt{h^{2s}-r^2}, d-1}]^2 \jump{u} \jump{v} (\xx)\dxx
\right|\\
&\le
\int_{f_0}
\left| B_{\sqrt{h^{2s}-r^2}, d-1} \jump{u}(\xx)
\frac{1}{2} \int_{B(\xx, \sqrt{h^{2s}-r^2})}
\nabla^2 \jump{v}(\xxi_\sss)|\sss-\xx|^2\dss\right|\\
&+                      
\left| 
B_{\sqrt{h^{2s}-r^2}, d-1} \jump{v}(\xx)
\frac{1}{2} \int_{B(\xx, \sqrt{h^{2s}-r^2})}
\nabla^2 \jump{u}(\xxi_\sss)|\sss-\xx|^2\dss\right|\\
&+
\left| 
\frac{1}{4} 
\int_{B(\xx, \sqrt{h^{2s}-r^2})}
\nabla^2 \jump{u}(\xxi_\sss)|\sss-\xx|^2\dss
\int_{B(\xx, \sqrt{h^{2s}-r^2})}
\nabla^2 \jump{v}(\xxi_\sss)|\sss-\xx|^2\dss\right|\dxx\\
&\lm
\max_{f} |\nabla^2 \jump{v}|
B_{\sqrt{h^{2s}-r^2},d-1}
\int_{f} |\jump{u}|(\xx)
\int_{B(\xx, \sqrt{h^{2s}-r^2})}|\sss-\xx|^2\dss\dxx\\
&+    
\max_{f} |\nabla^2 \jump{u}|
B_{\sqrt{h^{2s}-r^2},d-1}
\int_{f} |\jump{v}|(\xx)
\int_{B(\xx, \sqrt{h^{2s}-r^2})}|\sss-\xx|^2\dss\dxx\\
&+    
\max_{f} |\nabla^2 \jump{v}|
\max_{f} |\nabla^2 \jump{v}|
\int_{f_0}(h^{2s}-r^2)^{d+1}\dxx\\
&\lm
h^{-d-1}\int_{f} |\jump{v}| B_{\sqrt{h^{2s}-r^2},d-1}
        \int_{f} |\jump{u}|(h^{2s}-r^2)^{\frac{d+1}{2}}
+    
h^{-2d-2}\int_{f} |\jump{v}|
         \int_{f} |\jump{u}|(h^{2s}-r^2)^{d+1} \\
&\lm
h^{-d-1}(h^{2s}-r^2)^d 
\int_{f} |\jump{v}| \int_{f} |\jump{u}|
+
h^{-2d-2}(h^{2s}-r^2)^{d+1}
\int_{f} |\jump{v}| \int_{f} |\jump{u}|\\
&=
 (h^{-d-1}(h^{2s}-r^2)^d+ h^{-2d-2}(h^{2s}-r^2)^{d+1})
 \int_{f} |\jump{v}| \int_{f} |\jump{u}|.
\eali
\eees
In this way, we can estimate the last term in \eqref{4term_longest} as
\bees
\bali
&\frac{1}{[B_{h^s,d}]^2} 
\int_{-h^s}^{h^s}
\left|
\int_{f_0}
\int_{B(\xx,\sqrt{h^{2s}-r^2})} \jump{u} (\sss)\dss
\int_{B(\xx,\sqrt{h^{2s}-r^2})} \jump{v} (\sss)\dss 
-
[B_{\sqrt{h^{2s}-r^2},d-1}]^2 \jump{u} \jump{v} (\xx)\dxx
\right|\dr\\
&\lm
\frac{1}{[B_{h^s,d}]^2} 
\int_{-h^s}^{h^s}
(h^{-d-1}(h^{2s}-r^2)^d+ h^{-2d-2}(h^{2s}-r^2)^{d+1})
 \int_{f} |\jump{v}| \int_{f} |\jump{u}|\dr\\
&\lm
h^{-2sd} \int_{f} |\jump{v}| \int_{f} |\jump{u}|
\int_{-h^s}^{h^s}
h^{-d-1}(h^{2s}-r^2)^d+ h^{-2d-2}(h^{2s}-r^2)^{d+1})\dr\\
&\lm
h^{-2sd} \int_{f} |\jump{v}| \int_{f} |\jump{u}|
\cdot(h^{-d-1}h^{2sd+s} + h^{-2d-2}h^{2sd+3s})
\\
&=
(h^{-d-1+s} + h^{-2d-2+3s})
\int_{f} |\jump{v}| \int_{f} |\jump{u}|
\eali
\eees
and therefore, using \eqref{4term_longest} and the estimate 
\eqref{prop3_est} in Proposition \ref{prop2} we finally obtain
\bees
\bali
&\left|
(\eta_h*\jump{u}_f, \eta_h*\jump{v}_f) -
\frac{1}{[B_{h^s,d}]^2} \int_{-h^s}^{h^s} \int_{f}
[B_{\sqrt{h^{2s}-r^2},d-1}]^2 \jump{u}\jump{v} (\xx) 
\dxx\dr
\right|\\
&\lm
(h^{-d} + h^{-d-1+s} + h^{-2d-2+3s})
\int_{f} |\jump{v}| \int_{f} |\jump{u}|\\
&\lm
(h^{-d} + h^{-d-1+s} + h^{-2d-2+3s})h^d h^{s-1}
\|\nabla (\eta_h*u)\|_{\Omega} \|\nabla (\eta_h*v)\|_{\Omega}\\
&\lm
h^{s-1}(1+h^{-d+3s-2})
\|\nabla (\eta_h*u)\|_{\Omega} \|\nabla (\eta_h*v)\|_{\Omega}
\eali
\eees
as we have stated.\quad$\square$ 

\begin{cor}\label{cor_lemma7}
For all $u,v\in\mathbb{P}_{h,\kk}$ we have
$$
\bali
&\left|
\sum_{f\in\mcf} (\eta_h * \jump{u} ,\eta_h * \jump{v})_f -
\sum_{f\in\mcf} \frac{1}{[B_{h^s,d}]^2} \int_{-h^s}^{h^s} \int_{f}
[B_{\sqrt{h^{2s}-r^2},d-1}]^2 \jump{u}\jump{v} (\xx) 
\dxx\dr
\right|  \\
&\le
h^{s-1}(1+h^{3s-d-2}) \|\nabla(\eta_h*u)\| \|\nabla(\eta_h*v)\|.
\eali
$$
\end{cor}
Taking the sum of the inequalities in \eqref{est_last_aver} and applying the discrete Cauchy--Schwarz inequality 
$|\sum_{j\in J}a_j b_j|\le \sqrt {\sum_{j\in J}a_j} \sqrt{\sum_{j\in J}b_j}$ 
we obtain
$$
\bali
&
\left|
\sum_{f\in\mcf} (\eta_h * \jump{u} ,\eta_h * \jump{v})_f -
\sum_{f\in\mcf} \frac{1}{[B_{h^s,d}]^2} \int_{-h^s}^{h^s} \int_{f}
[B_{\sqrt{h^{2s}-r^2},d-1}]^2 \jump{u}\jump{v} (\xx) 
\dxx\dr
\right|\\
&\le
\sum_{f\in\mcf} \left|(\eta_h * \jump{u} ,\eta_h * \jump{v})_f -
\frac{1}{[B_{h^s,d}]^2} \int_{-h^s}^{h^s} \int_{f}
[B_{\sqrt{h^{2s}-r^2},d-1}]^2 \jump{u}\jump{v} (\xx) 
\dxx\dr
\right|\\
&\le
h^{s-1}(1+h^{3s-d-2}) 
\sum_{f\in\mcf}
\|\nabla(\eta_h*u)\|_{\Omega_+\cup\Omega_-} 
\|\nabla(\eta_h*v)\|_{\Omega_+\cup\Omega_-}\\
&\lm
h^{s-1}(1+h^{3s-d-2}) 
\|\nabla(\eta_h*u)\| \|\nabla(\eta_h*v)\|
\eali
$$
as stated in the corollary. \eproof

\emph{Remark:}
The above difference is lower order compared to  
$\|\nabla (\eta_h*u)\|\|\nabla (\eta_h*v)\|$
provided that $4s-d-3>0$ which is ensured for $s>1.5$.

Finally, we compute the approximation of the penalty term in 
\eqref{appr_last}, which appears in Lemma \ref{est_last_aver_lem}. 
\begin{itemize}
\item
For $d=2$ we have 
$$
\bali
&\frac{1}{[B_{h^s,2}]^2} \int_{-h^s}^{h^s} \int_{f}
[B_{\sqrt{h^{2s}-r^2},1}]^2 \jump{u}\jump{v} (\xx) 
\dxx\dr\\
&=
\frac{1}{h^{4s}\pi^2} \int_{-h^s}^{h^s} 4(h^{2s}-r^2)\dr 
                     \int_f  \jump{u}\jump{v} (\xx) \dxx =
\frac{16}{3\pi^2} h^{-s} \int_f  \jump{u}\jump{v}.
\eali
$$
\item
For $d=3$ we have 
$$
\bali
&\frac{1}{[B_{h^s,3}]^2} \int_{-h^s}^{h^s} \int_{f}
[B_{\sqrt{h^{2s}-r^2},2}]^2 \jump{u}\jump{v} (\xx) 
\dxx\dr\\
&=
\frac{9}{16 h^{6s}\pi^2} \int_{-h^s}^{h^s} \pi^2(h^{2s}-r^2)^2\dr 
                     \int_f  \jump{u}\jump{v} (\xx) \dxx 
=
\frac{3}{5} h^{-s} \int_f  \jump{u}\jump{v}.
\eali
$$
\end{itemize}

To prove the first main result we introduce the IP bilinear form
$a_{\IP,s}: \mathbb{P}_{h,\kk}\times \mathbb{P}_{h,\kk}\to\er$ 
with
\bee\label{IP_form_s}
a_{\textrm{IP,s}}(u,v) = 
(\nabla_h u, \nabla_h v) 
-\sum_{f\in\mathcal{F}}
(\aver{\nabla_h u}, \jump {v})_f + (\aver{\nabla_h v}, \jump {u})_f
+ \sum_{f\in\mathcal{F}}\sigma_{s,h} (\jump{u}, \jump {v})_f,
\eee
where
\bees
\sigma_{s,h} (\jump{u}, \jump {v})_f = 
\begin{cases}
\frac{16}{3\pi^2} h^{-s} (\jump{u}, \jump{v})_f \quad\textrm{for}\; d=2\\
\frac{3}{5} h^{-s} (\jump{u}, \jump{v})_f \quad\textrm{for}\; d=3.
\end{cases}
\eees
and the corresponding finite element approximation $u_{\IP,s}$ for which
\bee\label{IPS_approx}
a_{\IP,s}(u_{\IP,s}, v) = (g, v)\quad\forall v\in \mathbb{P}_{h,\kk}.
\eee 
\emph{Remark:}
Since we have the restriction $s>1.5$, the bilinear form $a_{\textrm{IP},s}$ 
can be recognized as an overpenalized IP bilinear form. 

\begin{thm}\label{thm1}
Assume that $3s>d+2$. Then the IP bilinear form in \eqref{IP_form_s} 
is a lower-order perturbation of $a_\eta$ in the sense that 
$$ 
|a_\eta(u,v) - a_{\IP,s} (u,v)| \lm 
h^{s-1}(1+h^{3s-d-2})
\|\nabla(\eta_h * u)\| \|\nabla(\eta_h * v)\|. 
$$
\end{thm}
\emph{Proof:}
Using Lemma \ref{diff_first_terms_lem}, Lemma \ref{23terms} 
and Corollary \ref{cor_lemma7} we obtain
$$
\bali 
&|a_\eta(u,v) - a_{\IP,s} (u,v)| \le
|(\eta_h*\nabla_h u, \eta_h*\nabla_h v) - (\nabla_h u, \nabla_h v)| \\
&+ 
\left|\sum_{f\in\mcf} 
(\eta_h*\eta_h*\nabla_h u, \jump{v})_f + (\eta_h*\eta_h*\nabla_h v, \jump{u})_f 
- (\aver{\nabla_h u}, \jump{v})_f - (\aver{\nabla_h v},  \jump{u}_f)\right|\\
&+
\left|
\sum_{f\in\mcf} (\eta_h * \jump{u} ,\eta_h * \jump{v})_f -
\sum_{f\in\mcf} \frac{1}{[B_{h^s,d}]^2} \int_{-h^s}^{h^s} \int_{f}
[B_{\sqrt{h^{2s}-r^2},d-1}]^2 \jump{u}\jump{v} (\xx) 
\dxx\dr
\right| \\
&\lm
(h^{s-1} + h^{s-1}(1+h^{3s-d-2}))
\|\nabla(\eta_h * u)\| \|\nabla(\eta_h * v)\|
\eali
$$
as stated in the theorem. \quad$\square$

Since the bilinear form $a_\eta$ is a slight modification of 
$a_{\IP,s}$ we expect that the local average of the approximations of 
$u_h$ and $u_{\IP,s}$ are also close to each other. In precise terms we 
have the following.
\begin{thm}\label{2ndthm}
Assume that $3s>d+2$. Then for the finite element approximations $u_h$ 
and $u_{\IP,s}$ we have
$$
\|\nabla(\eta_h*u_{\IP,s} - \eta_h*u_h)\| 
\lm 
h^{s-1} \|\nabla(\eta_h * u_h)\| + 
\max_j h_{\Omega_j}^d \|\eta_h * g - g_0\|.
$$
\end{thm}
\emph{Proof:}
Since $u_h$ solves \eqref{main} and $u_{\IP,s}\in \mathbb{P}_{h,\kk}$
we have 
$$
(\nabla(\eta_h * u_h), \nabla(\eta_h * (u_h - u_{\IP,s}))
=
(g_0, \eta_h * (u_h - u_{\IP,s}))
$$
such that using the equality 
$$
(\eta_h * w_1, w_2) =  (w_1, \eta_h * w_2) 
$$
for compactly supported functions $w_1, w_2\in L_1(\er^d)$ and the 
definition of $u_{\IP,s}$ in \eqref{IPS_approx} we obtain
\bee\label{long_for_thm2}
\bali
&
(\nabla(\eta_h * (u_h - u_{\IP,s})), \nabla(\eta_h * (u_h - u_{\IP,s})))\\
&=
(g_0, \eta_h * (u_h - u_{\IP,s})) - 
(\nabla(\eta_h * u_{\IP,s}), \nabla(\eta_h * (u_h - u_{\IP,s})))\\
&=
(g_0, \eta_h * (u_h - u_{\IP,s})) -
a_{\IP,s} (u_{\IP,s}, u_h - u_{\IP,s}) - 
(\nabla(\eta_h * u_{\IP,s}), \nabla(\eta_h * (u_h - u_{\IP,s})))\\ 
&+ a_{\IP,s} (u_{\IP,s}, u_h - u_{\IP,s})\\
&=
(\eta_h * g, u_h - u_{\IP,s}) - (g, u_h - u_{\IP,s}) - 
(\nabla(\eta_h * (u_{\IP,s} - u_h)), \nabla(\eta_h * (u_h - u_{\IP,s})))\\
&+ a_{\IP,s} (u_{\IP,s} - u_h, u_h - u_{\IP,s})
- 
(\nabla(\eta_h * u_h), \nabla(\eta_h * (u_h - u_{\IP})))
+ a_{\IP,s} (u_h, u_h - u_{\IP,s}). 
\end{aligned}
\eee
We note that the application of \eqref{normal_vs_conv_norm} to 
$u_h - u_{\IP,s}$ (instead of $\nabla_h u$) and the Friedrichs's 
inequality imply 
$$
\|u_h - u_{\IP}\|
\lm   
\|\eta_h*(u_h - u_{\IP})\|
\lm
\max_j h_{\Omega_j}^d
\|\nabla(\eta_h*(u_h - u_{\IP}))\|
$$
and therefore, using Theorem \ref{thm1} for the last two pair 
of terms in \eqref{long_for_thm2}, we obtain that
$$
\begin{aligned}
&
\|\nabla(\eta_h * (u_h - u_{\IP,s}))\|^2 \\
&\lm
\|\eta_h * g - g\| \|u_h - u_{\IP,s}\| + 
h^{s-1}(\|\nabla(\eta_h * (u_{\IP,s} - u_h))\|^2
+
\|\nabla(\eta_h * u_h)\|  \|\nabla(\eta_h * (u_h - u_{\IP,s}))\|)\\
&\lm
\max_j h_{\Omega_j}^d \|\eta_h * g - g\| 
\|\nabla(\eta_h * (u_h - u_{\IP,s}))\| + 
h^{s-1}(1+h^{3s-d-2})\|\nabla(\eta_h * (u_{\IP,s} - u_h))\|^2 \\
&+ 
h^{s-1}(1+h^{3s-d-2})
\|\nabla(\eta_h * u_h)\|  \|\nabla(\eta_h * (u_h - u_{\IP,s}))\|
\end{aligned}
$$
such that we finally get
$$
(1 - h^{s-1})\|\nabla(\eta_h * (u_h - u_{\IP,s}))\| 
\lm
\max_j h_{\Omega_j}^d \|\eta_h * g - g\| + 
h^{s-1} (1+h^{3s-d-2}) \|\nabla(\eta_h * u_h)\|,
$$
which implies the estimate in the theorem.\quad$\square$

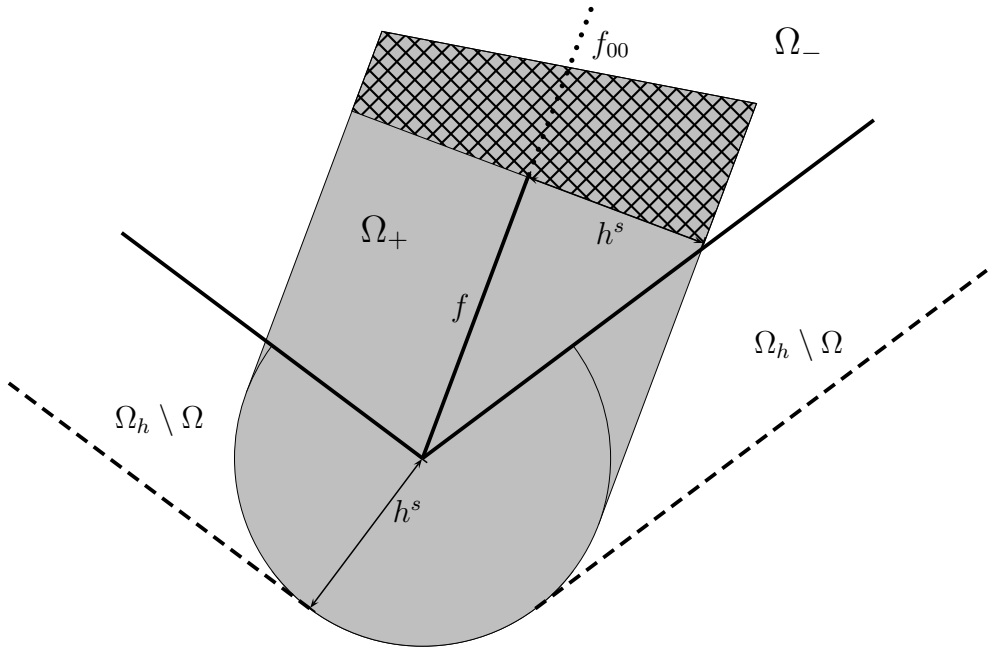
\begin{figure}
\begin{pspicture}(-1,0)(13,8) 

\pspolygon[linewidth=0cm,fillcolor=lightgray, fillstyle=solid](7.84,2.12)(9.94,7.72)(4.96,8.68)(3.16,3.88)
\pspolygon[linewidth=0cm, fillstyle=crosshatch](9.24,5.86)(9.94,7.72)(4.96,8.68)(4.56,7.62)
\psarc[linewidth=0cm,fillcolor=lightgray, fillstyle=solid](5.5,3){2.5}{143.4}{323.4}
\psarc[linewidth=0cm,fillcolor=lightgray, fillstyle=solid](5.5,3){2.5}{216.6}{36.6}


\psline[linewidth=0.05](5.5,3)(11.5,7.5)
\psline[linewidth=0.05](5.5,3)(1.5,6)
\psline[linewidth=0.05, linestyle=dashed](7,1)(13,5.5)
\psline[linewidth=0.05, linestyle=dashed](4,1)(0,4)

\psline[linewidth=0.05](5.5,3)(6.9, 6.74)
\psline[linewidth=0.07, linestyle=dotted](6.9, 6.74)(7.75,9)


\psline[linewidth=0.02]{|<->|}(4,1)(5.5,3)\rput(5.3,2.3){$h^s$}
\psline[linewidth=0.02]{|<->|}(6.9,6.74)(9.24,5.86)\rput(8,6){$h^s$}

\rput(8,8.5){$f_{00}$}
\rput(6,5){$f$}
\rput(5,6){\large $\Omega_+$}
\rput(10.5,8.5){\large $\Omega_-$}
\rput(2,3.5){$\Omega_h\setminus\Omega$}
\rput(10.5,4.5){$\Omega_h\setminus\Omega$}



\end{pspicture}
\caption{The support of $\eta_h*\jump{v_h}_f$
(shaded) and the set $\{f_{00}\otimes r: r\in [0, h^s]\}$ (with crosshatch)
in a 2-dimensional setup.}
\label{supports_jumps_out}
\end{figure}

To state quasi optimal convergence we observe that for each $h$ we have 
$\eta_h*u_h\in \mathbb{P}_{h,\kk,s}\subset H_0^1(\Omega_h)$. 
This means that the method in \eqref{main} is not conforming since the 
approximation in general is not in $H_0^1(\Omega)$. 
Also, the bilinear form $a_\eta^+$ is non-consistent in the sense that
the zero extension $u_0$ of $u$ is not necessarily the solution of \eqref{main}
for any $h$. In this way, we need to apply the Strang lemma \cite{ern04}, 
Section 2.3.2. For this we first note that for some constants $C_1$
and $C_2$ we have for all $0\not= w_1, w_2\in H_0^1(\Omega_h)$ that 
\bees
|a_\eta^+(w_1, w_2)| 
\le
C_1 \|w_1\|_{H^1_0(\Omega_h)} \|w_2\|_{H^1_0(\Omega_h)}
\eees
and 
\bees
C_2 \le \frac{a_\eta^+(w_1, w_2)}
{\|w_1\|_{H^1_0(\Omega_h)}\|w_2\|_{H^1_0(\Omega_h)}}. 
\eees
\begin{lem}\label{quasi_opt_lem}
The numerical solution $\eta_h * u_h$ of \eqref{main}
approximates $u$ in quasi optimal way in the sense that
\bees
\|\nabla(u- \eta_h*u_h)\|_{\Omega}\lm 
\inf_{v_h\in \mathbb{P}_{h,\kk}} \|\nabla(u- \eta_h*v_h)\|
+
h^{s-\frac{1}{2}} \|\nabla u\|.
\eees
\end{lem}
\emph{Proof:}
Since in this proof it is essential whether a scalar product 
is defined on $\Omega$ or on $\Omega_h$, we indicate it in the 
subscript.
A direct application of Lemma  2.25 in \cite{ern04} gives that
\bee\label{strang2_compare}
\bali
&\|\nabla(u- \eta_h*u_h)\|_{\Omega_h}\\
&\le 
(1+\frac{C_1}{C_2})
\inf_{v_h\in \mathbb{P}_{h,\kk}} \|\nabla(u- \eta_h*v_h)\|
+
\sup_{\eta_h*v_h\in \mathbb{P}_{h,\kk,s}} 
\frac{(\nabla u_0, \nabla(\eta_h*v_h))_{\Omega_h}
- (g_0, \eta_h*v_h)_{\Omega_h}}
{\|\eta_h*v_h\|_{1,\Omega_h}},
\eali
\eee
where the lower indices denote zero extensions. Using these,  
the second term in \eqref{strang2_compare} can be rewritten as
\bees
\bali
&
(\nabla u_0, \nabla(\eta_h*v_h))_{\Omega_h}-(g_0, \eta_h*v_h)_{\Omega_h}
=
(\nabla u, \nabla(\eta_h*v_h))_{\Omega}-(g, \eta_h*v_h)_{\Omega}\\
&=
(-\Delta u, \eta_h*v_h)_{\Omega}
+ \langle\nnu\cdot\nabla u, \eta_h*v_h\rangle_{\pa\Omega}
-(g, \eta_h*v_h)_{\Omega}
= \langle\nnu\cdot\nabla u, \eta_h*v_h\rangle_{\pa\Omega}.
\eali
\eees
Therefore, we can estimate \eqref{strang2_compare} to obtain
\bee\label{strang2_compare2}
\|\nabla(u- \eta_h*u_h)\|_{\Omega}\lm 
\inf_{v_h\in \mathbb{P}_{h,\kk}} \|\nabla(u- \eta_h*v_h)\|
+
\sup_{\eta_h*v_h\in \mathbb{P}_{h,\kk,s}} 
\frac{\langle\nnu\cdot\nabla u, \eta_h*v_h\rangle_{\pa\Omega}}
{\|\eta_h*v_h\|_{1,\Omega_h}}.
\eee
For the rest, it is sufficient to estimate the second term here. 
We apply a classical trace inequality in $\Omega$ and
in $\Omega_h\setminus \Omega$ which imply
\bee\label{otven}
\bali
&\sup_{\eta_h*v_h\in \mathbb{P}_{h,\kk,s}} 
\frac{\langle\nnu\cdot\nabla u, \eta_h*v_h\rangle_{\pa\Omega}}
{\|\eta_h*v_h\|_{1,\Omega_h}}\le
\sup_{\eta_h*v_h\in \mathbb{P}_{h,\kk,s}} 
\frac{\|\nnu\cdot\nabla u\|_{-\frac{1}{2},\pa\Omega} 
\|\eta_h*v_h\|_{\frac{1}{2},\pa\Omega}}
{\|\eta_h*v_h\|_{1,\Omega_h}}\\
&
\lm
\|u\|_{1,\Omega}
\sup_{\eta_h*v_h\in \mathbb{P}_{h,\kk,s}} 
\frac{\|\eta_h*v_h\|_{\frac{1}{2},\pa\Omega}}
{\|\eta_h*v_h\|_{1,\Omega_h}}
\lm
\|u\|_{1,\Omega}
\sup_{\eta_h*v_h\in \mathbb{P}_{h,\kk,s}} 
\frac{\|\eta_h*v_h\|_{1,\Omega_h\setminus\Omega}}
{\|\eta_h*v_h\|_{1,\Omega_h}}\\
&\lm
\|u\|_{1,\Omega}
\sup_{\eta_h*v_h\in \mathbb{P}_{h,\kk,s}} 
\frac{\|\nabla(\eta_h*v_h)\|_{\Omega_h\setminus\Omega}}
{\|\nabla(\eta_h*v_h)\|_{\Omega_h}}.
\eali
\eee
Observe that the numerator can be rewritten, using to Lemma 
\ref{lem1}, as 
$$
\|\nabla(\eta_h*v_h)\|_{\Omega_h\setminus\Omega} =
\|\eta_h*\nabla_h v_h + \eta_h*\jump{v_h}\|_{\Omega_h\setminus\Omega}
\le
\|\eta_h*\nabla_h v_h\|_{\Omega_h\setminus\Omega} 
+ \|\eta_h*\jump{v_h}\|_{\Omega_h\setminus\Omega}.
$$
To analyze the term $\|\eta_h*\jump{v_h}\|_{\Omega_h\setminus\Omega}$
we use the notations corresponding to Lemma \ref{compute_conv} 
and Fig. \ref{pict_ftrace}. Furthermore, we define 
$$
f_{00} = \{\xx\in f: d(\xx, \pa\Omega) > h^s\}
$$
such that 
$$
\textrm{supp}\: \eta_h*\jump{v_h}_f|_{\Omega_h\setminus\Omega}\subset
\textrm{supp}\: \eta_h*\jump{v_h}_f\setminus \{f_0\otimes r: r\in [0, h^s]\},
$$
see also Fig. \ref{supports_jumps_out}. 
The non-degeneracy of $\mct_h$ implies that the estimate in the 
remark after \eqref{fortriag1} is valid and therefore 
according to \eqref{fortriag1} the first term 
$\frac{1}{B_{h^s,d}^2}\int_{-h^s}^{h^s}I_1(r)\dr$ 
in the second line of \eqref{4term_longest} provides an upper bound 
for $\|\eta_h* \jump{v_h}_f\|_{\Omega_h\setminus\Omega}$. Therefore, 
the estimate in \eqref{1term_last} implies
$$
\|\eta_h* \jump{v_h}_f\|_{\Omega_h\setminus\Omega}^2
\le
h^{-d} \left[ \int_f \jump{v_h}\right]^2
\lm
h^{-d} h^d h^{s-1} \|\nabla (\eta_h*v_h)\|_{\Omega_+\cup\Omega_-}^2.
$$
Taking their sum for all interelement faces gives then 
\bee\label{strang_for1}
\|\eta_h* \jump{v_h}\|_{\Omega_h\setminus\Omega}
\lm
h^{\frac{s-1}{2}} \|\nabla (\eta_h*v_h)\|_{\Omega}.
\eee
Note that the condition on non-degeneracy implies that for all 
subdomains $\Omega_k\subset\tilde\Omega_j$ we have 
\bee\label{l_omegak}
\lambda(\Omega_k) \sim h_\Omega^d
\eee 
and 
\bee\label{l_omegajh}
\lambda(\Omega_{j,h}) \sim h^s h_\Omega^{d-1}.
\eee 
Using then \eqref{skala1}, \eqref{l_omegak}, \eqref{l_omegajh} and 
\eqref{normal_vs_conv_norm} we obtain that 
\bees
\bali
&\|\eta_h*\nabla_h v_h\|_{\Omega_{j,h}}^2 
\le 
\lambda(\Omega_{j,h}) \max_{\tilde\Omega_j} |\nabla_h v_h|^2
\lm
\lambda(\Omega_{j,h}) h_\Omega^{-d} |\nabla_h v_h|^2_{\tilde\Omega_j}\lm
h_\Omega^{s-1} |\nabla_h v_h|^2_{\tilde\Omega_j} \lm 
h_\Omega^{s-1} |\eta_h * \nabla v_h|^2_{\tilde\Omega_j}, 
\eali
\eees
which can be summed for all subdomain-patches to arrive at
\bee\label{strang_for2}
\|\eta_h*\nabla_h v_h\|^2_{\Omega_h\setminus\Omega} \lm
\sum_{\Omega_j\in\mathcal{T}_h} \|\eta_h*\nabla_h v_h\|_{\Omega_{j,h}}^2
\lm
h^{s-1} \sum_{\Omega_j\in\mathcal{T}_h} 
\|\eta_h*\nabla v_h\|_{\tilde\Omega_j}^2
\lm
h^{s-1}\|\eta_h*\nabla v_h\|^2.
\eee
We can use \eqref{strang_for1} and \eqref{strang_for2} to complete the 
estimation in \eqref{otven} as 
\bees
\sup_{\eta_h*v_h\in \mathbb{P}_{h,\kk,s}} 
\frac{\langle\nnu\cdot\nabla u, \eta_h*v_h\rangle_{\pa\Omega}}
{\|\eta_h*v_h\|_{1,\Omega_h}}\le
\|u\|_{1,\Omega}
\sup_{\eta_h*v_h\in \mathbb{P}_{h,\kk,s}} 
\frac{h^{s-\frac{1}{2}} \|\eta_h*\nabla v\|^2 }
{\|\nabla(\eta_h*v_h)\|_{\Omega_h}} 
\lm 
\|u\|_{1,\Omega}
h^{s-\frac{1}{2}}, 
\eees
which together with \eqref{strang2_compare2} gives the estimate in 
the lemma.
\eproof

We easily get now the statement on the convergence of the averaged 
IP method.

\begin{thm}\label{IP_quasiopt}
The averaged interior penalty approximation is quasi optimal
in the following sense:
\[
\|\nabla (u - \eta_h*u_{\IP,s})\| \lm 
 \inf_{v_h\in \mathbb{P}_{h,\emph{\kk}}} \|u - \eta_h*v_h\|_1 
+ \ordo (h^{s-\frac{1}{2}})
+ \max_j h_{\Omega_j}^d \|\eta_h*g - g_0\|.
\]
\end{thm}
\emph{Proof:}
A triangle inequality and the estimates 
in Theorem \ref{2ndthm} and Lemma \ref{quasi_opt_lem} imply that 
\[
\bali
\|\nabla (u - \eta_h*u_{\IP,s})\| &\lm 
\|\nabla (u - \eta_h * u_{h})\| +
\|\nabla(\eta_h * u_{\IP,s} - \eta_h * u_{h})\|\\
&\lm
 \inf_{v_h\in \mathbb{P}_{h,\kk}} \|u - \eta_h*v_h\|_1 +
\ordo(h^{s-\frac{1}{2}}) +  \max_j h_{\Omega_j}^d  \|\eta_h*g - g_0\|,
\eali
\]
as stated in the theorem. \quad$\square$




\emph{Remarks:}
The above derivation could cover the case of overpenalized IP methods
with $s>1.5$. The increase of the parameter $s$ can lead to 
ill-conditioned linear problems in the discretizations, such that
one should use appropriate preconditioners \cite{brenner08}.   

Based on the results of the paper, we propose the following introduction 
of IP methods for the numerical solution of \eqref{basic_eq}.
\begin{itemize}
\item
Introduce the $H^1$-conforming finite element discretization \eqref{main}.
\item
Since the $a_{\IP,s}$ bilinear form is a lower order approximation of 
$a_\eta$ and given more explicitly, one should compute $u_{\IP,s}$ in 
the practice.
\item
Compute the local average $\eta_h*u_{\IP,s}$. This converges to the weak 
solution 
$u$ of \eqref{basic_eq} in a quasi optimal way in the $H^1$-seminorm.
\end{itemize}

\section*{Appendix}
Following the notations in  \cite{evans10} we introduce
the smooth function $\Phi:\er^d\to\er$ with
$$
\Phi(\xx) =
\begin{cases}
C e^{\frac{1}{|\xx|^2-1}}\quad\textrm{if}\; |\xx|< 1\\
0 \quad\textrm{if}\; |\xx|\ge 1
\end{cases}
\quad\textrm{and}\quad 
\int_{B(\mathbf{0}, 1)} \Phi = 1
$$
and define $\Phi_\delta:\er^d\to\er$ by 
$\Phi_\delta := \left(\frac{1}{\delta}\right)^d \Phi(\frac{\xx}{\delta})$. 
Additionally, for $\xx\in\er^d$ we use the notation $\Phi_{\delta,\xx}$ 
for the function given simply by the  
$$
\Phi_{\delta,\xx} (\xx + \yy) :=  \Phi_{\delta} (\yy).
$$
We use the following proposition; for the proof we refer to
\cite{evans10}, pages 713--716.
\begin{prop}\label{appendix_prop}
For an arbitrary bounded Lipschitz domain $U$ and parameter $p\in [1,\infty)$
the following statements are valid.
\begin{itemize}
\item [(i)] 
For $f\in C(U)$ we have 
$\displaystyle{\lim_{\delta\to 0}\Phi_\delta * f\to f}$ uniformly.
\item [(ii)]
For any $f\in L_{p,\textrm{loc}}(U)$ we have 
$\displaystyle{\lim_{\delta\to 0}\Phi_\delta * f\to f}$ in 
$L_{p,\textrm{loc}}(U)$.
\end{itemize}
\end{prop}
We also need the following statements.
\begin{lem}\label{lem_tildef}
If for all $\xx\in K_1\cup K_2$ we have the limit
\bee\label{eq_for_tildef}
\lim_{\delta\to 0} (\eta_h*\jump{u}_f, \Phi_{\delta, \xx}) = \tilde f(\xx) 
\eee
then $\eta_h*\jump{u}_f$ can be identified with $\tilde f$. Also, 
for the function $\eta_h *_f \Phi_{\delta,\xx}:\mcf\to\er$ given by
$$
f\ni\yy\to\int_{\er^d} \eta_h(\yy-\zz)\Phi_{\delta,\xx}(\zz)\dzz
$$
we have the convergence 
\bee\label{convconv}
\lim_{\delta\to 0}\eta_h *_f \Phi_{\delta,\xx} 
= \eta_h (-\xx+\cdot)\quad
\textrm{in}\; L_1(f).
\eee
\end{lem}
\emph{Proof:}
We first note that 
$$
\bali
&\langle \eta_h*\jump{u}_f, \Phi_{\delta, \xx}\rangle
=
\int_{\er^d} \eta_h*\jump{u}_f(\xx+\yy) \Phi_{\delta, \xx}(\xx+\yy)\: 
\mathrm{d}\yy
= 
\int_{\er^d} \eta_h*\jump{u}_f(\xx+\yy) \Phi_{\delta}(\yy)\:\mathrm{d}\yy\\
&=
\int_{\er^d} \eta_h*\jump{u}_f(\xx+\yy) \Phi_{\delta}(-\yy)\:\mathrm{d}\yy
=
\Phi_\delta *\eta_h*\jump{u}_f (\xx).
\eali
$$
In this way, according to property (ii) we can rewrite 
the condition in \eqref{eq_for_tildef} as 
\bee\label{assu2}
\lim_{\delta\to 0} \Phi_\delta *\eta_h*\jump{u}_f \to \tilde f 
\quad\textrm{in}\;L_{1,\textrm{loc}}(\er^d). 
\eee
Using the property $(i)$ above, the fact that  $\eta_h*\jump{u}_f$ is 
locally integrable and the limit in \eqref{assu2},
we have that for each function $g\in C_0^\infty(\Omega)$ the 
following equality is valid:
\[
\bali
\langle\eta_h*\jump{u}_f, g\rangle 
&= (\eta_h*\jump{u}_f, g)
 = \lim_{\delta\to 0} (\eta_h*\jump{u}_f, \Phi_\delta *g)
 = \lim_{\delta\to 0} (\Phi_\delta *\eta_h*\jump{u}_f, g)\\
 = (\tilde f, g), 
\eali
\]
which proves the first statement of the lemma.

To prove the second statement we rewrite 
$\eta_h *_f \Phi_{\delta,\xx}$ as 
$$
\eta_h *_f \Phi_{\delta,\xx} (\yy) 
=  
\int_{\er^d} \eta_h(\yy-\zz) \Phi_\delta (-\xx+\zz)\dzz
=
\int_{\er^d} \eta_h(\yy-\xx-\zz) \Phi_\delta (\zz)\dzz.
$$
Accordingly, we have the pointwise convergence 
$$
\eta_h *_f \Phi_{\delta,\xx} (\yy) \to  \eta_h(\yy-\xx).
$$
On the other hand, 
$$
\int_{\er^d} \eta_h(\yy-\xx-\zz) \Phi_\delta (\zz)\le 
\max_{\er^d} |\eta_h|
$$
so that the function $\max_{\er^d}|\eta_h|\cdot\mathbbm{1}\in L_1(\mcf)$ 
delivers an upper bound for each function $\eta_h *_f \Phi_{\delta,\xx}$.
The statement is therefore an obvious consequence of the Lebesgue
dominant convergence theorem.
\quad$\square$\\\\
\emph{Proof of Lemma \ref{compute_conv}:}
We compute $\eta_h*\jump{u}_f$ based on the first statement in Lemma 
\ref{lem_tildef}.
For this, we use Lemma \ref{lem1} and \eqref{convconv}, which give 
$$
\begin{aligned}
& \lim_{\delta\to 0} \langle\eta_h*\jump{u}_f, \Phi_{\delta,\xx}\rangle
=\lim_{\delta\to 0} \langle\jump{u}_f, \eta_h*\Phi_{\delta,\xx}\rangle
=\lim_{\delta\to 0} \int_f\jump{u}_f(\yy) \eta_h*_f\Phi_{\delta,x}(\yy)\dyy\\
&=\int_{f} \jump{u}_f(\yy) 
\lim_{\delta\to 0} \eta_h*_f\Phi_{\delta,x}(\yy)\dyy 
=
\int_{f} \jump{u}_f(\yy)  \eta_h(\yy-\xx)\dyy
\end{aligned}
$$
as stated in the lemma. \quad $\square$

\emph{Proof of Proposition \ref{prop2}}
We first prove \eqref{prop2_est}.
According to Lemma \ref{lem1} and the consecutive remark, 
we have obviously that 
\bee\label{prop2_1}
\int_{f_0} |\jump{v}| \le |v|_\BV,
\eee
where the BV seminorm is taken on $\overline{K_{+}\cup K_{-}}$. 
For the next step we use a scaling argument and introduce the 
function space 
$$
\bar{\mathbb{P}}_{K} := \mathbb{P}_{h,\kk}|_{K_{+}\cup K_{-}} \diagup 
\langle\mathbbm{1}\rangle,
$$  
which is the restriction of $\mathbb{P}_{h,\kk}$ to $K_{+}\cup K_{-}$ 
factorized with the constant functions. The BV seminorm on 
this function space becomes a norm, and accordingly, we use 
the notation $\|\cdot\|_{\BV}$. We next prove that for all
$\epsilon>0$ there is $h_0>0$ such that for all $h<h_0$ and 
$v\in\bar{\mathbb{P}}_{K}$ we have
\bee\label{second_step_ineq}
\|\eta_h*v - v\|_\BV < \epsilon \|v\|_\BV.
\eee
For this we consider a normed basis $\{v_1, v_2, \dots, v_D\}$ 
of $v\in\bar{\mathbb{P}}_{K}$ with respect to the $\BV$ norm 
and define the Euclidean norm $\|\cdot\|_E$
generated by this basis such that 
\bee\label{eucl_norm}
\left|\sum_{j=1}^D a_j v_j \right|_E ^2 = \sum_{j=1}^D a_j ^2.
\eee
This norm should be equivalent with the $\BV$ norm, i.e. there 
is a constant $c_0$ with 
\bee\label{norm_eq}
\|v\|_E \le c_0 \|v\|_\BV \quad\forall\; v\in\bar{\mathbb{P}}_{K_0}.
\eee
Note that this constant should be not the same for all pair
of neighboring subdomains, but it is a continuous function of the 
position of the vertices. In particular, if we fix the edge $f_0$
of length \emph{one}, then $f_0$ is fixed for $d=2$ and for $d=3$ the 
remaining vertex should be in a compact set depicted in fig ...
if the condition of non-degeneracy holds true.
Therefore, the constant $c_0$ has a finite maximum. Similarly, if 
we fix now an arbitrary interelement face chosen above the remaining 
node of $K_{0-}$ and $K_{0+}$ can lie in a compact set. In this way,
for each pair of neighboring subdomains with at least one interelement
edge of length \emph{one} there is a uniform constant $c_0$ in 
\eqref{norm_eq}.   
Also, since we have a finite basis, and $\eta_h$ is a Dirac 
series, there is $h_0$ such that for all $h<h_0$ we have 
\bee\label{conv_small}
\|\eta_h * v_j - v_j\|_\BV \le \frac{\epsilon}{c_0D} 
\|v_j\|_\BV = \frac{\epsilon}{c_0\sqrt{D}}  \quad\forall\; j\in\{1,2,\dots,D\}.
\eee
We obtain also here that \eqref{conv_small} is valid for all
pair of neighboring subdomains with at least one interelement
edge of length \emph{one} with a uniform parameter $h_0$.
Then using \eqref{conv_small}, \eqref{eucl_norm} and \eqref{norm_eq}
we have that for any $0<h<h_0$ and 
$v= \sum_{j=1}^D a_j v_j \in \bar{\mathbb{P}}_{K}$
the following inequality is valid:
$$
\bali
&\|\eta_h*v -v\|_{\BV} = \|\sum_{j=1}^D \eta_h*v_j -v_j\|_{\BV}\le
\sum_{j=1}^D \|\eta_h*a_jv_j -a_jv_j\|_{\BV}\\ 
&
\le \frac{\epsilon}{c_0\sqrt{D}} \sum_{j=1}^D |a_j|
\le \frac{\epsilon}{c_0} \sqrt{\sum_{j=1}^D |a_j|^2} =
 \frac{\epsilon}{c_0} \|v\|_E \le \epsilon \|v\|_\BV, 
\eali
$$
which proves the inequality in \eqref{second_step_ineq}.
Consequently, we also have 
\bees
(1-\epsilon) \|v\|_\BV \le \|v - \eta_h * v\|_\BV +  
\|\eta_h * v\|_\BV - \epsilon \|v\|_\BV \le \|\eta_h * v\|_\BV.
\eees  
In the last step we relate \eqref{prop2_1} and 
\eqref{second_step_ineq} and use that $\eta_h * v$ 
is differentiable to obtain
\bee\label{pre_lem}
\int_{f_0} \jump{v} \le \|v\|_\BV \le \frac{1}{1-\epsilon}\|\eta_h * v\|_\BV 
= \frac{1}{1-\epsilon} \int_{K_{+}\cup K_{-}} |\nabla (\eta_h* v)|.
\eee
To prove the statement of the lemma for two arbitrary neighboring 
subdomains $\Omega_+$ and $\Omega_-$ we use \eqref{pre_lem} and the 
equalities in \eqref{basic_scale2} which give  
\bee\label{pre_lem_final}
\bali
&\int_{f_\Omega} \jump{v} = h_\Omega^{d-1} \int_{f_0} \jump{v_0} 
\lm h_\Omega^{d-1} \int_{K_{+}\cup K_{-}} |\nabla (\eta_{h_0}* v_0)|\\
&
= h_\Omega^{-d} h_\Omega^{d-1} \int_{\Omega_{+}\cup\Omega_-} h_\Omega^{1} 
|\nabla (\eta_{h_0h_\Omega^{\frac{1}{s}}} * v)|
= \int_{\Omega_{+}\cup \Omega_-} 
|\nabla (\eta_{h_0 h_\Omega^{\frac{1}{s}}} * v)|. 
\eali
\eee
The inequality remains true if the lower index $h_0 h_\Omega^{\frac{1}{s}}$
is changed to a smaller one since this is equivalent with the choice of a
smaller index $h_0$. Obviously the condition $h^{1-\frac{1}{s}}<h_0$ implies  
$$
h_0 h_\Omega^{\frac{1}{s}} > h
$$
and using  \eqref{pre_lem_final} with the previous remark gives that 
$$
\int_{f_\Omega} \jump{v} 
\le
\int_{\Omega_{+}\cup \Omega_-} 
|\nabla (\eta_{h} * v)| 
$$
as stated in the inequality \eqref{prop2_est}. 

To prove \eqref{prop3_est} we again use the geometric setup discussed 
at the beginning of Section \ref{sect_4}, where for the simplicity, 
we use the notation $h=h_\Omega$. 
With this we obtain
\bee\label{est_aver_2norm}
\bali
&\|\eta_h*\nabla u\|_{h\cdot K_+\cup h\cdot K_-}
\left( \int_{h\cdot f_0 \dot{\pm} h^s} 1\right)^{\frac{1}{2}}
\ge 
\|\eta_h*\nabla u\|_{h\cdot f_0 \dot{\pm} h^s} 
\left( \int_{h\cdot f_0 \dot{\pm} h^s} 1\right)^{\frac{1}{2}}\\
&=
\|\nabla (\eta_h* u)\|_{h\cdot f_0 \dot{\pm} h^s}
\left( \int_{h\cdot f_0 \dot{\pm} h^s} 1\right)^{\frac{1}{2}}
\ge
\int_{h\cdot f_0 \dot{\pm} h^s} |\nabla (\eta_h* u)|\\
&\gm
\int_{h\cdot f^*}\left|
\int_{-h^s}^{h^s}\nabla (\eta_h* u)(x,\yy)\dx\right|\dyy
=
\int_{h\cdot f^*}\left|
(\eta_h* u)(h^s,\yy) - (\eta_h* u)(-h^s,\yy)
\right|\dyy\\
&\ge
 \int_{h\cdot f^*}\left|u(h^s,\yy) - u(-h^s,\yy) \right| 
- \left|(\eta_h* u)(h^s,\yy) - u(h^s,\yy) \right|
- \left|u(-h^s,\yy) - (\eta_h* u)(-h^s,\yy) \right|\dyy.
\eali
\eee
To continue with the estimate we note that $u$ is differentiable twice 
in $B((-h^s,\yy),h^s)$ and according to \eqref{midp_aver} we have
$$
\bali
&
\left|(\eta_h* u)(h^s,\yy) - u(h^s,\yy)\dyy \right|\le
\frac{1}{B_{h^s,d}}\cdot\frac{1}{2}\max_{h\cdot K_-}|\nabla^2u|
\int_{B(0,h^s)}|\sss|^2\dss\\
&\lm
\max_{h\cdot K_-}|\nabla^2u| h^{-sd}h^{s(d+2)} = 
\max_{h\cdot K_-}|\nabla^2u| h^{2s}
\eali
$$
Therefore, using \eqref{est_aver_2norm} we have
\bees
\|\eta_h*\nabla u\|_{h\cdot K_+\cup h\cdot K_-}
\left( \int_{h\cdot f_0 \dot{\pm} h^s} 1\right)^{\frac{1}{2}}
\gm
\int_{h\cdot f^*}\left|\jump{u} \right|
- h^{2s} \lambda(h\cdot f^*) \max_{h\cdot K_+\cup h\cdot K_-}|\nabla^2u|,
\eees 
which can be rewritten with the aid of \eqref{max_nabla_sq_2_norm} 
and the condition $s\ge \frac{3}{2}$ as
$$
\bali
&\int_{h\cdot f^*}\left|\jump{u} \right|
\lm
h^{2s} h^{d-1} \max_{h\cdot K_+\cup h\cdot K_-}|\nabla^2u|
+
h^{\frac{s}{2}+\frac{d}{2}-\frac{1}{2}}
\|\eta_h*\nabla u\|_{h\cdot K_+\cup h\cdot K_-}\\
&\lm
h^{2s} h^{d-1} h^{-\frac{d}{2}-1} \|\nabla u\|_{h\cdot K_+\cup h\cdot K_-}
+
h^{\frac{s}{2}+\frac{d}{2}-\frac{1}{2}}
\|\eta_h*\nabla u\|_{h\cdot K_+\cup h\cdot K_-}\\
&\lm
h^{\frac{d}{2}}(h^{\frac{s}{2}-\frac{1}{2}}+h^{2s-2})
\|\eta_h*\nabla u\|_{h\cdot K_+\cup h\cdot K_-}
\lm
h^{\frac{d}{2}}h^{\frac{s}{2}-\frac{1}{2}}
\|\eta_h*\nabla u\|_{h\cdot K_+\cup h\cdot K_-}
\eali
$$
A summation with respect to the faces gives then the desired inequality.
\quad $\square$

\section*{Acknowledgments}
The author acknowledges the support of the Hungarian Research Fund OTKA 
(grants PD10441 and K104666).

\bibliography{/neumann/u/staff/applanal/izsakf/ferenc_bib.bib}

\def\cprime{$'$}
\begin{thebibliography}{10}

\bibitem{arnold01}
D.~N. Arnold, F.~Brezzi, B.~Cockburn, and L.~D. Marini.
\newblock Unified analysis of discontinuous {G}alerkin methods for elliptic
  problems.
\newblock {\em SIAM J. Numer. Anal.}, 39(5):1749--1779 (electronic), 2001/02.

\bibitem{attouch06}
H.~Attouch, G.~Buttazzo, and G.~Michaille.
\newblock {\em Variational analysis in {S}obolev and {BV} spaces}, volume~6 of
  {\em MPS/SIAM Series on Optimization}.
\newblock Society for Industrial and Applied Mathematics (SIAM), Philadelphia,
  PA, 2006.
\newblock Applications to PDEs and optimization.

\bibitem{brenner08}
S.~C. Brenner, L.~Owens, and L.-Y. Sung.
\newblock A weakly over-penalized symmetric interior penalty method.
\newblock {\em ETNA. Electronic Transactions on Numerical Analysis [electronic
  only]}, 30:107--127, 2008.

\bibitem{buffa09}
A.~Buffa and C.~Ortner.
\newblock Compact embeddings of broken {S}obolev spaces and applications.
\newblock {\em IMA J. Numer. Anal.}, 29(4):827--855, 2009.

\bibitem{cockburn03b}
B.~Cockburn, M.~Luskin, C.-W. Shu, and E.~S{\"u}li.
\newblock Enhanced accuracy by post-processing for finite element methods for
  hyperbolic equations.
\newblock {\em Math. Comp.}, 72(242):577--606 (electronic), 2003.

\bibitem{cockburn98}
B.~Cockburn and C.-W. Shu.
\newblock The local discontinuous {G}alerkin method for time-dependent
  convection-diffusion systems.
\newblock {\em SIAM J. Numer. Anal.}, 35(6):2440--2463 (electronic), 1998.

\bibitem{csorgo14}
G.~Cs\"org\H{o} and F.~Izs\'ak.
\newblock Energy norm error estimates for averaged discontinuous {G}alerkin
  methods in 1 dimension.
\newblock {\em Int. J. Numer. Anal. Model.}, 2014.
\newblock to appear.

\bibitem{dawson11}
C.~Dawson, E.~J. Kubatko, J.~J. Westerink, C.~Trahan, C.~Mirabito, C.~Michoski,
  and N.~Panda.
\newblock Discontinuous {G}alerkin methods for modeling hurricane storm surge.
\newblock {\em Advances in Water Resources}, 34(9):1165--1176, 2011.

\bibitem{dipietro10}
D.~A. Di~Pietro and A.~Ern.
\newblock Discrete functional analysis tools for discontinuous {G}alerkin
  methods with application to the incompressible {N}avier-{S}tokes equations.
\newblock {\em Math. Comp.}, 79(271):1303--1330, 2010.

\bibitem{dipietro12}
D.~A. Di~Pietro and A.~Ern.
\newblock {\em Mathematical Aspects of Discontinuous Galerkin Methods}.
\newblock Springer-Verlag, Berlin, Heidelberg, 2012.

\bibitem{ern04}
A.~Ern and J.-L. Guermond.
\newblock {\em Theory and {P}ractice of {F}inite {E}lements}, volume 159 of
  {\em Applied Mathematical Sciences}.
\newblock Springer-Verlag, New York, 2004.

\bibitem{evans10}
L.~C. Evans.
\newblock {\em Partial differential equations}, volume~19 of {\em Graduate
  Studies in Mathematics}.
\newblock American Mathematical Society, Providence, RI, second edition, 2010.

\bibitem{gudi10}
T.~Gudi.
\newblock A new error analysis for discontinuous finite element methods for
  linear elliptic problems.
\newblock {\em Math. Comp.}, 79(272):2169--2189, 2010.

\bibitem{halmos50}
P.~R. Halmos.
\newblock {\em Measure {T}heory}.
\newblock D. Van Nostrand Company, Inc., New York, N. Y., 1950.

\bibitem{hartmann10}
R.~Hartmann, J.~Held, T.~Leicht, and F.~Prill.
\newblock Discontinuous galerkin methods for computational aerodynamics. 3d
  adaptive flow simulation with the {D}{L}{R} {P}{A}{D}{G}{E} code.
\newblock {\em Aerospace Science and Technology}, 14(7):512 -- 519, 2010.

\bibitem{hesthaven07}
J.~Hesthaven and T.~Warburton.
\newblock {\em Nodal Discontinuous Galerkin Methods: Algorithms, Analysis, and
  Applications}, volume~54 of {\em Texts in Applied Mathematics}.
\newblock Springer, New York, first edition, 2007.

\bibitem{hirsch99}
F.~Hirsch and G.~Lacombe.
\newblock {\em Elements of Functional Analysis}, volume 192 of {\em Graduate
  Texts in Mathematics}.
\newblock Springer-Verlag, New York, 1999.
\newblock Translated from the 1997 French original by Silvio Levy.

\bibitem{houston07b}
P.~Houston, D.~Sch\"otzau, and T.~P. Wihler.
\newblock Energy norm a posteriori error estimation of hp-adaptive
  discontinuous {G}alerkin methods for elliptic problems.
\newblock {\em Mathematical Models and Methods in Applied Sciences},
  17(01):33--62, 2007.

\bibitem{ji12}
L.~Ji, Y.~Xu, and J.~K. Ryan.
\newblock Accuracy-enhancement of discontinuous {G}alerkin solutions for
  convection-diffusion equations in multiple-dimensions.
\newblock {\em Math. Comp.}, 81(280):1929--1950 (electronic), 2012.

\bibitem{karakashian03}
O.~A. Karakashian and F.~Pascal.
\newblock A posteriori error estimates for a discontinuous {G}alerkin
  approximation of second-order elliptic problems.
\newblock {\em SIAM J. Numerical Analysis}, 41(6):2374--2399, 2003.

\bibitem{king12}
J.~King, H.~Mirzaee, J.~K. Ryan, and R.~M. Kirby.
\newblock Smoothness-{I}ncreasing {A}ccuracy-{C}onserving ({SIAC}) {F}iltering
  for {D}iscontinuous {G}alerkin {S}olutions: Improved {E}rrors {V}ersus
  {H}igher-{O}rder {A}ccuracy.
\newblock {\em J. Sci. Comput.}, 53:129--149, 2012.

\bibitem{mirzaee13a}
H.~Mirzaee, J.~King, J.~Ryan, and R.~Kirby.
\newblock Smoothness-increasing accuracy-conserving filters for discontinuous
  {G}alerkin solutions over unstructured triangular meshes.
\newblock {\em SIAM Journal on Scientific Computing}, 35(1):A212--A230, 2013.

\bibitem{mirzaee14}
H.~Mirzaee, J.~K. Ryan, and R.~M. Kirby.
\newblock Smoothness-increasing accuracy-conserving ({SIAC}) filters for
  discontinuous {G}alerkin solutions: Application to structured tetrahedral
  meshes.
\newblock {\em Journal of Scientific Computing}, 58(3):690--704, 2014.

\bibitem{peletier07}
M.~A. Peletier, R.~Planqu{\'e}, and M.~R{\"o}ger.
\newblock Sobolev regularity via the convergence rate of convolutions and
  {J}ensen's inequality.
\newblock {\em Ann. Sc. Norm. Super. Pisa Cl. Sci. (5)}, 6(4):499--510, 2007.

\bibitem{riviere08}
B.~Rivi{\`e}re.
\newblock {\em Discontinuous {G}alerkin methods for solving elliptic and
  parabolic equations}, volume~35 of {\em Frontiers in Applied Mathematics}.
\newblock Society for Industrial and Applied Mathematics (SIAM), Philadelphia,
  PA, 2008.
\newblock Theory and implementation.

\bibitem{szucs13}
Z.~Sz{\"u}cs.
\newblock The {L}ebesgue decomposition of representable forms over algebras.
\newblock {\em J. Operator Theory}, 70(1):3--31, 2013.

\bibitem{tago12}
J.~Tago, V.~M. Cruz-Atienza, J.~Virieux, V.~Etienne, and F.~J. S\'anchez-Sesma.
\newblock A 3{D} $hp$-adaptive discontinuous {G}alerkin method for modeling
  earthquake dynamics.
\newblock {\em Journal of Geophysical Research: Solid Earth}, 117(B9):n/a--n/a,
  2012.

\end{thebibliography}

\bibliographystyle{abbrv}

\end{document}